\documentclass[12pt,a4paper]{article}

\usepackage{amsmath}
\usepackage{amssymb}

 \newcommand{\bs}[1]{\boldsymbol{#1}}
 \newcommand{\floor}[1]{\lfloor{#1}\rfloor}
 \newcommand{\ceiling}[1]{\lceil{#1}\rceil}

 \newcommand{\RR}{\mathbb R}
 
 \newcommand{\ZZ}{\mathbb Z}

 \newcommand{\cM}{\mathcal M}
 \newcommand{\cP}{\mathcal P}
 \newcommand{\cB}{\mathcal B}

 \newcommand{\es}{\varnothing}

 \newcommand{\s}{{\scriptscriptstyle <}}
 \newcommand{\g}{{\scriptscriptstyle >}}

 \newtheorem{theorem}{Theorem}
 \newtheorem{lemma}{Lemma}
 \newtheorem{prop}{Proposition}
 \newtheorem{coro}{Corollary}
 \newtheorem{fact}{Fact}

 \newcommand{\qed}{\hfill $\square$}
 \newcommand{\btimes}{\mbox{\huge \raisebox{-0.2ex}{$\times$}}}

\begin{document}

\centerline{\Large \bf An exactly solved model for}

\bigskip
\centerline{\Large \bf mutation, recombination and selection}

\vspace{15mm}
        
\centerline{{\sc Michael Baake}
        and {\sc Ellen Baake}}
\vspace{5mm} 

{\small
\begin{center}
    Institut f\"ur Mathematik und Informatik,\\
    Universit\"at Greifswald, Jahnstr.~15a,\\
    17487 Greifswald, Germany  
\end{center}
}

\vspace{10mm}
\begin{abstract} 
It is well known that rather general mutation-recombination models can
be solved algorithmically (though not in closed form) by means of Haldane 
linearization. The price to be paid is that one has to work with a
multiple tensor product of the state space one started from.

Here, we present a relevant subclass of such models, in continuous time,
with independent mutation events at the sites, and crossover
events between them. It admits a closed solution of the
corresponding differential equation on the basis of the
original state space, and also closed expressions for the linkage 
disequilibria, derived by means of M\"obius inversion.
As an extra benefit, the approach can be extended to a model
with selection of additive type across sites.
We also derive a necessary and sufficient criterion for the mean fitness 
to be a Lyapunov function and determine the asymptotic behaviour of the 
solutions.
\end{abstract}

\vspace{10mm}
\hphantom{Key Words:}           
measure-valued dynamical systems, M\"obius inversion


\clearpage

\section*{Introduction}

The basic mechanisms which create genetic variation in biological
evolution are {\em mutation\/} and {\em recombination\/}.
They are counteracted by {\em selection\/}, which removes variation. 
Genetic information may be quite generally described in terms of
a collection of linearly ordered  sites (i.e.\ a sequence of sites), each
of which is occupied by an element of a given (finite or infinite)
set which we denote as {\em site space\/}; if this set is finite,
it is often termed {\em alphabet\/}. A specific sequence is also
called {\em type\/}.

{\em Mutation\/} is  treated as a random state change of a site variable, 
which occurs independently at every site.
\!\!{\em Recombination\/} occurs on the occasion of sexual reproduction, and
refers to the creation of `offspring' sequences from two 
(randomly chosen) `parental' ones, where a subset of the `maternal'
sites is combined with the complementary set of the `paternal'
sites, and the linear ordering along the sequence is maintained.
This process is realized through one, or a number of,
{\em crossover events\/}, where the two parental strands are
interlaced between a pair of neighbouring sites.
An important feature of recombination is that it removes dependencies
between sites, known as {\em linkage disequilibria\/} in genetics.
Finally, {\em selection\/} is caused by the flourishing of fit
individuals at the expense of less fit ones.

We consider an infinite population of sequences which evolves
under the joint action of mutation, selection or recombination,
or of any combination thereof.
This is to be considered as the infinite population limit (IPL) 
of the stochastic process alluded to, and defines a deterministic
dynamical system for probability measures
(in discrete or in continuous time). It describes the time
evolution of the measure with probabilistic certainty, see
\cite[Ch.\ 11]{EK}, and Thm.\ 2.1 of it in particular. 
Although there are many interesting and
important questions connected with finite populations, we 
focus on the differential equation of the deterministic limit here,
which we will call {\em IPL equation\/} from now on.
In particular, we will not employ the traditional discrete dynamical
systems, but follow the continuous route along the lines of
Kimura \cite{Kimura} and Akin \cite{Akin-Buch}, which happens 
to be much less developed than
it ought to be, see also \cite{BG} for a recent review.

Mutation is a linear process and straightforward to
deal with. Selection involves some nonlinearity, 
which is due to norm conservation under the dynamics,
but this nonlinearity may be removed through a simple transformation.
Recombination contains a very different source of nonlinearity,
which is due to the fact that pairs of objects  are
involved in the process, and is much harder to treat.
Nevertheless, if both state space and time are
discrete, a procedure (known as {\em Haldane linearization\/},
see \cite{MHR,Dawson,Dawson2} and \cite[Ch.\ 6]{Lyu}) 
is available which transforms the dynamical system
(exactly) into a linear one. It involves a multilinear transformation
of the probabilities to a new set of variables, namely
certain linkage disequilibria,
which describe the deviation from statistical independence
of  sites. These variables  decay independently
and geometrically, i.e.\ they decouple and diagonalize
the dynamics. 
Unfortunately, however, the procedure is cumbersome since it
relies on recursions, and no closed form is
available for the transformation in the general case.

In a previous paper \cite{EB}, the special case of single crossover 
events was considered, where offspring sequences are composed
from one maternal and one paternal segment. 
This scenario is particularly relevant in molecular evolution,
where crossover events are rare, and it is most consistently described
in continuous time. For discrete site spaces, and with the help of
the corresponding vector space structure, the
linearizing transform could be given in closed form with the
help of elementary methods from multilinear algebra.

The aim of this article is to further develop this approach in a 
systematic measure-theoretic setting
which also incorporates more general site spaces and does not
require an explicit coordinatization. We will essentially start
from the deterministic IPL equation and construct its solution
explicitly, first for recombination only. The so-called
M\"obius inversion principle will then give a simple approach to
the calculation of a suitable (and, in particular,
complete) set of linkage disequilibria.
It will then turn out that mutation and even selection may be included in
the framework, provided fitness is additive, meaning that
the fitness of any type may be decomposed into a sum of
independent contributions of its individual sites, i.e.\ if
there is no interaction between sites. Such  results may 
be helpful for the solution of the corresponding inverse problem,
i.e.\ the determination of recombination rates from experimental
data, e.g.\ observed patterns of linkage disequilibria
along sequences \cite{Clark,Schaeffer}.

The exposition will be more explicit than needed for a purely mathematical 
audience, and we also try to give rather precise references to 
background material we use. We hope that the article will
become more self-contained this way and that it is also accessible
for readers with a more biological background. 

The structure of
the paper is as follows. After some preliminaries in Section \ref{prel},
we will briefly summarize the description of mutation through an IPL
equation on the space of positive measures in Section \ref{mut},
followed by some general remarks on measure-valued IPL equations.
The core of the article is Section \ref{recombi}, where we solve, step 
by step, the IPL equation for recombination and construct an explicit 
solution of the abstract Cauchy problem, together with a closed form of 
the corresponding linkage disequilibria. The latter is based on an
application of the inclusion-exclusion principle via M\"obius inversion
(a supplement is given in the Appendix).
Section \ref{both} combines mutation and recombination. Section
\ref{sel} deals with selection and recombination, with some emphasis on
the role of mean fitness as a Lyapunov function. Finally, Section \ref{all}
ties together all three evolutionary forces --- still giving an explicit 
solution, expressions for the linkage disequilibria, and asymptotic
properties. We close with some  afterthoughts
mainly aimed at  the relationship to models in discrete time.

\bigskip
\bigskip
\section{Preliminaries} \label{prel}

If $X$ is a locally compact space (by which we always mean to include the 
Hausdorff property), we use $\cM_+(X)$ to denote the set of 
{\em finite\/} positive regular Borel measures on $X$, with
$0\in\cM_+(X)$. Likewise, $\cM(X)$ is the vector space of real (or signed) 
finite regular Borel measures. It is a Banach space under the norm 
$\|\omega\| = |\omega|(X)$ where $|\omega|$ denotes total variation 
measure. Due to the Riesz-Markov representation theorem, $\cM(X)$ can also
be viewed as the dual of $C_{\infty}(X,\RR)$, the Banach space of real-valued
continuous functions which vanish at infinity, equipped with the
usual supremum norm, see \cite[Thm.~IV.18]{RS}, as well as
\cite[Ch.\ 6]{Rudin} and \cite[Ch.\ IV.4]{RS} for general 
background material. Note that $\cM(X)$ with the variation norm 
$\|.\|$ is actually a Banach lattice, and this gives access to the highly 
developed theory of positive operators \cite{Sch,Arendt}. We will mainly 
be interested in the closed convex subsets 
$\cM_+^m(X) := \{\omega\in\cM^{}_+(X) \mid \omega(X)=m\}$,
and in $\cP(X)=\cM_+^1(X)$ in particular, the set of {\em probability
measures\/} on $X$. Note that, for positive measures $\omega$, we simply
have $\|\omega\| = \omega(X)$.

If the Borel $\sigma$-algebra of $X$ is generated by a family of sets
that is closed under finite intersections, a regular Borel measure on $X$
is already uniquely specified by its values on the elements of this generating
family \cite{Berb,Lang}. This is a property that we will need several times,
in particular if $X=X_1\times X_2$ is a product space, equipped with
the product topology. 
\begin{fact} \label{baire}
   Let $\nu,\nu'$ be two regular Borel measures on the locally compact
   product space $X=X_1\times X_2$ which coincide on all ``rectangles'' 
   $E_1\times E_2$ where $E_1$ and $E_2$ each run through 
   the Borel sets of $X_1$ and $X_2$. Then $\nu=\nu'$, i.e.\
   $\nu(E)=\nu'(E)$ for all Borel sets\/ $E$ of\/ $X$.
\end{fact}

{\sc Proof}: In view of the above remark, the only obstacle to cope with is
the (non-vacuous!) situation when the $\sigma$-algebra generated by the 
rectangles $E_1\times E_2$ is {\em not\/} the full Borel $\sigma$-algebra 
of $X$. However, the $\sigma$-algebra 
generated by the rectangles contains all Baire sets $F$ of $X$, because the
Baire sets of $X$ possess the required Cartesian product property
\cite[Lemma 56.2]{Berb}, and the Borel sets of $X_i$ contain the Baire sets
of $X_i$. The equality of $\nu$ and $\nu'$ now follows from 
\cite[Thm.\ 62.1]{Berb} (this rests on the fact that every Baire measure
has a unique extension to a regular Borel measure).  \qed

\smallskip
Standard examples of locally compact spaces include the compact ones, such 
as any finite set or the closed interval $[0,1]$, but also $\RR^k$ and 
$\ZZ^{\ell}$ with $k,\ell\ge 0$, and arbitrary combinations thereof. 
These certainly cover all meaningful parameter spaces to be expected in 
biological applications.

If $X$ is a finite set (which is an important case in population
genetics), $\cP(X)$ is a {\em simplex\/}. If the cardinality
of $X$ is $M$, this simplex has dimension $M-1$, i.e.\ any probability
measure is a unique convex linear combination of the $M$ {\em extremal\/}
measures that constitute the vertices of the simplex. If
$X=\{1,\dots , M\}$, they are denoted by $e^{}_i$, $i=1,\dots ,M$,
and fixed by their values on singleton sets, 
$e^{}_i(\{j\}) = \delta_{i,j}$. In other words, any $\omega\in\cP(X)$
is of the form $\omega = \sum_{i=1}^{M} a^{}_i e^{}_i$
with all $a_i\ge 0$ and $a^{}_1 + \ldots + a^{}_M = 1$.
This provides the canonical coordinatization of this situation.

The set (or state space) $X$ that we need will have a product structure, 
described on the basis of {\em sites\/}. For later convenience, we use
$N=\{0,1,\dots ,n\}$ for the set of sites, i.e.\ we start counting
with $0$ here. To site $i$, we attach the locally 
compact space $X_i$, and our state space is then
\begin{equation} \label{seq-space}
   X \; = \; X_0 \times X_1 \times \ldots \times X_n\,,
\end{equation}
which is still locally compact. One Banach space of measures to show 
up is the space $\cM(X)$ with the corresponding variation norm $\|.\|$.
Note that $\cM(X)$ contains the (algebraic) tensor product space
$\cM^{\otimes}:=\bigotimes_{i=0}^n \cM(X_i)$, and also its completion
(here, the closure in the given $\|.\|$-norm of $\cM(X)$).
To simplify notation, the latter will also be denoted by $\cM^{\otimes}$,
because we shall only deal with Banach spaces here. Recall that 
$\cM^{\otimes}$ contains the product measures 
$\omega=\omega^{}_0\otimes\dots\otimes\omega^{}_n$ with
$\omega^{}_i\in\cM(X_i)$, but also all (finite) linear combinations
of measures of this kind. Because we consider the completion, also
all measures are contained which can be approximated with such
linear combinations in the norm. All probability measures
of product form are in this space, but note that the single measures
in the product need not be probability measures themselves.

If $X_i=\{1,\dots ,M_i\}$ is finite, for all $0\le i\le n$, $X$
is still a finite set, with $M=\prod_{i=0}^n M_i$ elements.
Then, $\cM(X) = \cM^{\otimes}$, and this is simply 
a real vector space of dimension $M$. $\cM(X) = \cM^{\otimes}$
is also true for $X$ discrete. In this case, the action of
operators in tensor product form is well defined. In general, 
if $\cM^{\otimes}\subsetneq\cM(X)$, one can still go beyond 
$\cM^{\otimes}$ under certain circumstances, e.g.\ by including
integrals (rather than finite sums) of product measures. However,
we do not want to enter this rather technical discussion, and refer
to \cite[Ch.~IX.6]{Lang} and \cite[Ch.~IV.7]{Sch} for some
background material, and to \cite[Ch.~13]{Dudley} for some of
the problems that are related to these difficulties.

$X$ finite is the case most frequently studied in the theory of sequence
evolution, and it was the motivation for this work, see \cite{EB}
and references given there. However, many results hold in greater
generality, which we want to cover in view of potential applications in 
quantitative genetics. There, the space $X_i$ often is a state
space such as $\RR$, or a compact subset thereof. In this case,
$\cM^{\otimes}$ is a true subspace of $\cM(X)$, which
has to be taken care of later on (occasional restrictions of $X$ to a 
finite set will be mentioned explicitly).

The main reason for using the above set $N$ of sites is that we will need
{\em ordered partitions\/} of $N$, which are uniquely specified by a
set of cuts or {\em crossovers\/}. The possible cut positions are at 
the {\em links\/} between sites, which we denote by half-integers, i.e.\ 
by elements of the
set $L=\{\frac{1}{2}, \frac{3}{2}, \dots  ,\frac{2n-1}{2}\}$.
We will use Latin indices for sites and Greek indices for links,
and the implicit rule will always be that $\alpha=\frac{2i+1}{2}$
is the link between site $i$ and $i+1$.

With this notation, the ordered partitions of $N$ are in one-to-one
correspondence with the subsets of $L$ as follows. If 
$A=\{\alpha^{}_1,\dots ,\alpha^{}_p\}\subset L$, let $N^{}_A$ denote 
the ordered partition
$$\{0,\dots ,\floor{\alpha^{}_1}\} \, , \,
  \{\ceiling{\alpha^{}_1},\dots ,\floor{\alpha^{}_2}\}\, , \; \dots  \; ,
  \{\ceiling{\alpha^{}_{p}},\dots , n \} $$
where $\floor{\alpha}$ ($\ceiling{\alpha}$) denotes the largest
integer below $\alpha$ (the smallest above $\alpha$). In particular,
we have $N^{}_{\es} = N$ and $N^{}_L = \{\{0\},\dots,\{n\}\}$.
With this definition, it is clear that $N^{}_B$ is a {\em refinement\/} 
of $N^{}_A$ if and only if $A\subset B$. Consequently, the lattice of 
ordered partitions of $N$  corresponds to the Boolean algebra of the
finite set $L$, denoted by $\cB(L)$, cf.\ \cite[Ch.\ I.2]{Aigner}. We
prefer this notation to that with partitions, as it is easier to 
deal with. If $A\subset B$, we will write $B\!-\!A$ for $B\setminus A$, 
and $\stackrel{\;\_\!\_}{A}$ for the set $L\!-\!A$.

This setup allows us to use the powerful tool of {\em M\"obius inversion\/}
from combinatorial theory \cite[Ch.\ IV.2]{Aigner}, which is a systematic way 
to employ the inclusion-exclusion principle. If $f$ and $g$ are mappings from 
$\cB(L)$ to $\RR$ which are, for all $A\subset L$, related by
\begin{equation} \label{moebius1}
   g(A) \; = \; \sum_{B\subset A} f(B) \, ,
\end{equation}
then this can be solved for $f$ via the inversion formula 
\cite[Thm.\ 4.18]{Aigner}
\begin{equation} \label{moebius2}
   f(A) \; = \; \sum_{B\subset A} g(B)\, \mu(B,A)
\end{equation}
with the M\"obius function $\mu(B,A) = (-1)^{|A-B|}$, where $|A\!-\!B|$
stands for the cardinality of the set $A\!-\!B$. For $B$ not a subset
of $A$, we set $\mu(B,A) = 0$ which makes the M\"obius function
into an element of the so-called incidence algebra, see 
\cite[Ch.\ IV.1]{Aigner} for details. It is important to note that M\"obius 
inversion is not restricted to functions, it also applies to bounded operators.

\bigskip
\bigskip
\section{Mutation and Markov generator} \label{mut}

The description of mutation is rather straight-forward.
Let us start from a {\em finite\/} population. 
Since we are working in continuous time, we assume an
independent Poisson clock for each individual member of a (finite)
population, and a mutation occurs for an individual whenever
its clock rings, according to prescribed mutation rates between
(finitely many) types or states. Since the individuals are
independent, this is a simple Markov process for each of them.
If we now go to the infinite population limit, the time evolution 
of the probability measure for the types is, almost surely, 
described by a (deterministic) ordinary differential equation (ODE).
This is the so-called {\em IPL equation\/}, compare
\cite[Thm.\ 11.2.1]{EK} for a general justification, which we will
also rely on below. For the simple
mutation case, this ODE is linear. It clearly coincides with the
ODE for the probability measure of the individual Markov process,
usually obtained from multiple realizations through the 
law of large numbers.

Let us consider the case that $X$ is a finite state space of 
cardinality $|X|=M$ in more detail, where
$\dim^{}_{\RR}(\cM(X))=M$. The mutation rate from state $\ell$ to state 
$k$ is given by $Q^{}_{k,\ell} = Q^{}_{k\leftarrow\ell}$, where we
already consider $Q$ as a mapping acting on the corresponding
probabilities, resp.\ measures. The rate matrix $Q$ is a Markov generator,
i.e.\ it has non-negative entries everywhere except on its diagonal, and 
vanishing column\footnote{In contrast to the standard probability literature,
we adopt the transposed version here since we are 
considering the situation from the (linear) operator point of view.}
sums. The time evolution is then 
fully described by the Markov semigroup $\{\exp(tQ)\mid t\ge 0\}$, see
\cite[Ch.\ 1.1 and Ch.\ 4.2]{EK}. We shall usually assume that $Q$ is 
{\em irreducible\/}, 
i.e.\ it is possible to reach every state from any other one. In this case, the 
equilibrium state is unique and given by the properly normalized 0-eigenvector 
of the generator $Q$. It can actually be given in closed form,
see \cite[Lemma 6.3.1]{FW}.

If $X$ has the product structure introduced above, our mutation process is
supposed to be of a more special form, for biological reasons. We assume that
mutation happens at all sites in parallel and independently from one another,
so that our generator has the form
\begin{equation} \label{markov1}
   Q \; = \; \sum_{i=0}^{n} \, Q^{}_i
\end{equation}
where each $Q^{}_i$ is, in a properly coordinatized way, the tensor product of
a rate matrix at site $i$ and unit matrices of matching dimension everywhere 
else, i.e.\
\begin{equation} \label{qi}
 Q^{}_i \; = \; \bs{1}^{}_{M_0}\otimes\dots\otimes\bs{1}^{}_{M_{i-1}} \otimes 
   q^{}_i \otimes \bs{1}^{}_{M_{i+1}}\otimes\dots\otimes\bs{1}^{}_{M_n} 
\end{equation}
where $q^{}_i$ is a local rate matrix (of dimension $M_i$) for the state space 
$X_i$, acting on $\cM(X_i)$.
The rate matrices $Q^{}_i$ clearly commute with one another. Note also that 
$Q$ of (\ref{markov1}) is irreducible if and only if all the $q^{}_i$ are.
The Markov semigroup inherits the tensor product structure,
i.e.\ we have
\begin{equation} \label{tensor1}
\exp(tQ) \; = \; \prod_{i=0}^{n} \exp(tQ^{}_i)
         \; = \; \bigotimes_{i=0}^{n} \exp(tq^{}_i)\, .
\end{equation}

In view of our following description of recombination, we prefer to avoid
an explicit coordinatization here, so we will not use matrix notation.
This simply means that we have to reinterpret the
generator $Q$ as a linear operator on $\cM(X)$. Nothing of the above
actually changes, we only have to read $Q$ (or $q^{}_i$) as a linear mapping
on $\cM(X)$ (or on $\cM(X_i)$). The two conditions for $Q$ to be a Markov 
generator now read as follows (the analogous conditions apply to $q^{}_i$ in 
relation to $\cM(X_i)$).
\begin{enumerate}
\item If $\nu$ is a positive measure and
     $E$ any Borel set such that $\nu(E)=0$, then $(Q\nu)(E)\ge 0$.
\item If $\nu$ is a positive measure, then
     $(Q\nu)(X) = 0$.
\end{enumerate}
The first condition ensures that the semigroup generated by $Q$ maps 
$\cM_+(X)$ into itself. Under the present circumstances, where $Q$ is
bounded and $\cM(X)$ is a reflexive Banach space, this condition is 
necessary and sufficient for the positivity of $\exp(tQ)$, see 
\cite[Thm.\ 1.11]{Arendt}. It is sometimes also called
the {\em positive minimum principle\/}. 
The second condition means that the semigroup is Markov, i.e.\ it preserves 
the norm of positive measures, and, in particular, maps $\cP(X)$ into 
itself. In this setting, irreducibility implies that the 
kernel of the Markov generator $Q$ is one-dimensional.

The IPL equation for our simple mutation process\footnote{In this linear
case, the IPL equation is closely related to the master equation
commonly used in the physics literature, see \cite[Ch.~5]{vankamp} for 
details.} now reads
\begin{equation} \label{m-master}
   \dot{\omega} \; = \; \Phi^{}_{\rm mut}(\omega)
   \; := \; \Bigl(\, \sum_{i=0}^{n} \, Q^{}_i \Bigr) \,\omega
\end{equation}
which we will take, in generalization of the discrete situation, as the starting
point for the analysis of mutation, without tracing it back to an explicit
stochastic process. We then obtain, by employing standard results \cite{Amann} 
from the theory of ordinary linear differential equations in (finite-dimensional)
Banach spaces (see also Theorem \ref{general-ODE} below):
\begin{prop} \label{markov-sol}
   The abstract Cauchy problem of the IPL equation $(\ref{m-master})$ 
   with initial condition\/ $\omega^{}_0\in\cP(X)$ has the unique solution
$$ \omega^{}_t \; = \; \exp\Bigl(t\sum_{i=0}^{n}\, Q^{}_i \Bigr)\,\omega^{}_0 $$
   which is, for $t\ge 0$, a one-parameter family of probability measures. \qed
\end{prop}

To formulate a generalization of Prop.~\ref{markov-sol},
let us forget about the product structure for a moment and consider
the linear ODE  $$\dot{\omega} \; = \; Q \omega$$ with $Q$ the generator 
of a uniformly (or norm) continuous Markov semigroup on $\cM(X)$,
compare \cite[Ch.\ I.3]{Engel}. 
This is the case if and only if the linear operator $Q$, in addition to
satisfying assumptions 1.\ and 2.\ from above, is bounded, and hence
defined on all of $\cM(X)$, see \cite[Cor.\ II.1.5]{Engel}.
In particular, we can then write the semigroup in exponential
form \cite[Thm.\ I.3.7]{Engel}, i.e.\ as $\exp(tQ)$, and the solution 
as $\,\omega^{}_t = \exp(tQ)\, \omega^{}_0$. 
In what follows, we will (non-constructively) assume that a process is given 
that leads to a bounded generator $Q$ which is a linear operator on $\cM(X)$,
i.e.\ maps regular Borel measures to
regular Borel measures. As long as this is the case, it is sufficient
to work with assumptions 1.\ and 2., even if the space of measures
considered is no longer reflexive. The analogue of Prop.~\ref{markov-sol}
then holds on the Banach subspace $\cM^{\otimes}$, to which we shall
restrict our attention whenever $Q$ is of the form specified in
Eqs.\ (\ref{markov1}) and (\ref{qi}). This makes no difference at
all as long as $X$ is discrete.

Many results can still be generalized to densely defined generators of 
strongly continuous semigroups, see \cite[Ch.\ I.5]{Engel}, but already 
the well-posed\-ness of the Cauchy problem needs some thought, compare
\cite[Ch.\ II.6]{Engel} for a discussion. Also, the characterization of
generators for positive semigroups becomes more involved, see 
\cite[Ch.\ 3]{Arendt}. Usually, one would then rather describe the entire
process by means of semigroups on function spaces, compare \cite[Ch.\ 1.4]{EK}.
Since all explicit mutation schemes we have in mind lead to uniformly continuous 
semigroups, we will not expand on the more general situation.

Let us instead add a few remarks on the general type of IPL equation 
that arises when recombination and selection are also included. 
This will also better explain our formulation of mutation,
from the point of view of measure-valued differential equations.
In what follows, it is sufficient to investigate the first order ODE
\begin{equation} \label{gen-master}
   \dot{\omega} \; = \; \Phi(\omega)
\end{equation}
on the Banach space $\cM(X)$, where $\Phi$ is a mapping from $\cM(X)$ into 
itself (alternatively, we can study (\ref{gen-master}) on any closed
subspace of $\cM(X)$ that is invariant under $\Phi$). 
Unlike $\Phi^{}_{\rm mut}$ from (\ref{m-master}), $\Phi$ need not be 
linear, and it is the nonlinear cases below that we are most interested in. 
The three properties we will meet below are:
\begin{itemize}
\item[{\bf A1}] The mapping $\Phi$ is (globally) Lipschitz.
\item[{\bf A2}] If $\nu\in\cM_+(X)$, i.e.\ $\nu$ is a positive measure,
   and $E$ any Borel set such that $\nu(E)=0$, then we have 
   $\big(\Phi(\nu)\big)(E) \ge 0$.
\item[{\bf A3}] For any $\nu\in\cM_+(X)$, we have
   $\big(\Phi(\nu)\big)(X) = 0$.
\end{itemize}
It is clear that our formulation of mutation constitutes a linear example
of such a mapping. 

\begin{theorem} \label{general-ODE}
   If\/ $\Phi\!: \cM(X)\to\cM(X)$ satisfies\/ {\rm (A1)}, the abstract Cauchy
   problem of the\/ {\rm ODE} $(\ref{gen-master})$, with initial condition\/ 
   $\omega^{}_0\in\cM(X)$, has a unique solution. 
   If\/ $\Phi$ also satisfies\/ {\rm (A2)}, the cone $\cM_+(X)$ 
   of positive measures is invariant under the semiflow for\/ $t\ge 0$
   $($in other words, $\cM_+(X)$ is positive invariant\/$)$.
   {}Finally, if\/ $\Phi$ also satisfies\/ {\rm (A3)}, the norm of positive
   measures is preserved in forward time. In particular, the convex set\/
   $\cP(X)$ of probability measures is then positive invariant.
\end{theorem}
{\sc Proof}: If $\Phi$ is Lipschitz, we can invoke the Picard-Lindel\"of
Theorem for ODEs on Banach spaces, see \cite[Thm.\ 7.6]{Amann}, so existence
and uniqueness of the solution of the abstract Cauchy problem are clear.

If $\Phi$ also satisfies (A2), positive invariance of $\cM_+(X)$ follows 
from a continuity argument, see p.~235 and Thm.~16.5 together with Remark~16.6 
of \cite{Amann} for a proof. If $\Phi$ is linear, (A2) is the so-called
positive minimum principle, and our assertion also follows from
\cite[Thm.\ 1.11]{Arendt}, which uses a functional analytic proof.

{}Finally, assume $\Phi$ satisfies (A1) -- (A3).
Let $\omega^{}_0\in\cM_+^m(X)$ be the initial condition and denote the 
corresponding unique solution of (\ref{gen-master}) by $\omega^{}_t$.
Then, $\omega^{}_t\in\cM_+(X)$ for all $t\ge 0$ by the previous argument, 
so $\|\omega^{}_t\|=\omega^{}_t(X)$. This implies 
$\frac{{\rm d}}{{\rm d}t}\|\omega^{}_t\| = \big(\Phi(\omega^{}_t)\big)(X)
= 0$ by assumption (A3), so $\|\omega^{}_t\|\equiv\|\omega^{}_0\|=m$.
This proves the assertion. \qed

\bigskip
\bigskip
\section{Recombination} \label{recombi}

This section deals with the nonlinear IPL equation for recombination,
and is the core of our article.  We develop the 
results step by step here. The combination with mutation will then be rather 
painless, and an addition of selection will be discussed after that.

\subsection{Recombination on measures}

Let $X,Y$ be two locally compact spaces with attached 
measure spaces $\cM(X)$ and $\cM(Y)$. If $f\!: X\to Y$ is a continuous 
function and $\omega\in\cM(X)$, then $f.\omega := \omega\circ f^{-1}$
is an element of $\cM(Y)$, where $f^{-1}(y) :=
\{x\in X\mid f(x)=y\}$ means the preimage of $y\in Y$ in $X$, with obvious
extension to $f^{-1}(B)$, the preimage of a subset $B\subset Y$ in $X$.
Due to the continuity of $f$, $f^{-1}(B)$ is a Borel set in $X$
if $B$ is a Borel set in $Y$.

Let $X=X_0\times\ldots\times X_n$ be as in Section \ref{prel},
and let, from now on, $N$ and $L$ always denote the set of sites and 
links as introduced there. In this section, we can entirely work
with the Banach space $\cM(X)$, equipped with the variation norm $\|.\|$.
Let $\pi^{}_i \! : X \to X_i$ be the canonical
{\em projection\/} which is continuous. It induces a mapping from $\cM(X)$ 
to $\cM(X_i)$ by $\omega\mapsto\pi^{}_i .\omega$, where
$(\pi^{}_i .\omega) (E) = \omega(\pi^{-1}_i (E))$,
for any Borel set $E\subset X_i$. By (slight) abuse
of notation, we will use the symbol $\pi^{}_i$ also for this induced mapping.
It is clear that $\pi_i$ is linear and maps positive measures to positive
measures of the same norm. As such, it is bounded and hence also continuous.
In particular, it maps $\cP(X)$ to $\cP(X_i)$
and may then be understood as marginalization. Likewise, we can
start from any (ordered) index set $I\subset N$ and define a projector
$\pi^{}_I\!: \cM(X) \to \cM(X_I)$ with $X_I:=\btimes_{\!i\in I}\,X_i$.
With this notation, $X_N=X$. We will frequently also use the
abbreviation $\pi^{}_{\s\alpha}$ for the projector 
$\pi^{}_{\{1,... ,\floor{\alpha}\}}$, and $\pi^{}_{\g\alpha}$ 
for $\pi^{}_{\{\ceiling{\alpha},... ,n\}}$. These objects
may be understood as `cut and forget' operators, since they give
the distribution of what is left after a cut is made at $\alpha$, and the
trailing resp.\ leading segment is discarded.

This now enables us to introduce the elementary {\em recombination
operator\/}, or {\em recombinator\/} as we will call it from now on,
$R_{\alpha}\! : \cM(X)\to\cM(X)$, for $\alpha\in L$.
If $\omega = 0$, $R_{\alpha} (\omega) := 0$, and otherwise
\begin{equation} \label{reco-op1}
   R_{\alpha} (\omega) \; := \; \frac{1}{\|\omega\|}\,
     \bigl((\pi^{}_{\s\alpha}.\omega) \otimes 
           (\pi^{}_{\g\alpha}.\omega)\bigr)
\end{equation}
which is a (partial) product measure. Here and in what follows, we tacitly
identify (if necessary) a product measure with its unique extension to a
regular Borel measure on $X$, which is justified by Fact \ref{baire}.
The following property is now an immediate consequence of the definition.
\begin{fact} \label{positivity}
   The recombinator $R_{\alpha}$ maps $\cM_+(X)$ into itself and 
   preserves the norm of positive measures. In particular,
   it maps\/ $\cP(X)$ into itself. \qed
\end{fact}

Let us comment on the choice of (\ref{reco-op1}). 
Being composed of the cut-and-forget operators
for the leading and the trailing ends, $R_{\alpha} (\omega)$
has the interpetation of a `cut-and-relink operator', which
describes a cut at $\alpha$, followed by (random) reunion
of the resulting segments.  

At first sight,
it might appear more natural to drop the prefactor $1/\|\omega\|$. 
However, the norm of a positive measure $\omega$ would then not be 
preserved unless $\|\omega\|=1$. In view of later extensions, it
is more desirable not to be restricted to probability
measures, and that is why we prefer (\ref{reco-op1}) which makes
$R_{\alpha}$ positive homogeneous of degree 1,
\begin{equation} \label{pos-hom}
   R_{\alpha} (a\,\omega) \; = \; |a| \, R_{\alpha}(\omega) \, ,
\end{equation}
for arbitrary $a\in\RR$. Note, however, that $R_{\alpha}$ is {\em not\/} 
a linear operator, not even when restricted to $\cM_+(X)$.

\begin{fact} \label{lipschitz}
   Let $\alpha\in L$. The recombinator $R_{\alpha}$ satisfies
   $\|R_{\alpha}(\omega)\|\le\|\omega\|$, for all $\omega\in\cM(X)$,
   and is $($globally\/$)$ Lipschitz on $\cM(X)$.
\end{fact}
{\sc Proof}: 
Let us first observe that, for arbitrary $\omega,\omega'\in\cM(X)$ and
$\alpha\in L$, we obtain the inequality
$$ \|(\pi^{}_{\s\alpha}.\omega)\otimes(\pi^{}_{\g\alpha}.\omega')\|
   \; \le \; \|\omega\|\,\|\omega'\|\,,$$
which is a simple consequence of Hahn's decomposition for real measures,
see \cite[Thm.\ 6.14]{Rudin}, applied separately to the factors of the 
product measure. For $0\neq \omega\in\cM(X)$, we then have
$$ \|R_{\alpha}(\omega)\| \; = \;
   \frac{\|(\pi^{}_{\s\alpha}.\omega)\otimes(\pi^{}_{\g\alpha}.\omega)\|}
        {\|\omega\|} \; \le \; \|\omega\| \, , $$
with equality for positive measures, as stated in Fact \ref{positivity}.
Clearly, we also have $R_{\alpha}(0) = 0$, so that the first assertion
follows.

Let $\omega,\omega'\in\cM(X)$. If one of them is the 0-measure, say
$\omega' = 0$, we have 
$ \|R_{\alpha}(\omega) - R_{\alpha}(\omega')\|
= \|R_{\alpha}(\omega)\|\le\|\omega\|=\|\omega -\omega'\|$. 
So we may assume both $\omega$ and $\omega'$ to be
different from 0 and hence to have positive norm.
With the above inequalities, we can now employ the following
3$\varepsilon$-type argument
\begin{eqnarray*}
  \lefteqn{
  \| R_{\alpha}(\omega) - R_{\alpha}(\omega')\| } \qquad \\[3mm]
  & = &
  \bigg\|\, \frac{\bigl(\pi^{}_{\s\alpha}.\omega\bigr)
  \otimes\bigl(\pi^{}_{\g\alpha}.(\omega-\omega')\bigr)}{\|\omega\|}
   + \, \frac{\bigl(\pi^{}_{\s\alpha}.(\omega-\omega')\bigr)
  \otimes\bigl(\pi^{}_{\g\alpha}.\omega'\bigr)}{\|\omega'\|} \\[2mm]
  & & \quad + \, \left(\frac{1}{\|\omega\|} - \frac{1}{\|\omega'\|}\right)\,
  \Bigl[\bigl(\pi^{}_{\s\alpha}.\omega\bigr)\otimes
        \bigl(\pi^{}_{\g\alpha}.\omega'\bigr)\Bigr] \, \bigg\|  \\[3mm]
  & \le &  \frac{\big\|\bigl(\pi^{}_{\s\alpha}.\omega\bigr)
  \otimes\bigl(\pi^{}_{\g\alpha}.(\omega-\omega')\bigr)\big\|} 
  {\|\omega\|} \, + \,
  \frac{\big\|\bigl(\pi^{}_{\s\alpha}.(\omega-\omega')\bigr)
  \otimes\bigl(\pi^{}_{\g\alpha}.\omega'\bigr)\big\|}{\|\omega'\|} \\[2mm]
  & & \quad + \,
  \bigg| \frac{1}{\|\omega\|} - \frac{1}{\|\omega'\|}\bigg|\cdot
  \big\|\bigl(\pi^{}_{\s\alpha}.\omega\bigr)\otimes
  \bigl(\pi^{}_{\g\alpha}.\omega'\bigr)\big\| \\[2mm]
  & \le & 2 \, \| \omega -\omega'\| + 
  \big| \|\omega\| - \|\omega'\|\big|
  \;\; \le \;\; 3\, \|\omega - \omega'\|\, .
\end{eqnarray*}
Together, this gives the second assertion,
with Lipschitz constant $\le 3$. \qed

\smallskip
In view of Fact \ref{positivity}, it makes sense to investigate the
properties of the recombinators restricted to the positive cone
$\cM_+(X)$. The crucial property which
underlies our later analysis is the following.
\begin{prop} \label{comm-rules1}
   The elementary recombinators, when restricted to $\cM_+(X)$, are 
   idempotents and commute with one another. In other words, we then have
   $R_{\alpha}^2=R_{\alpha}$ and
   $R_{\alpha}R_{\beta}=R_{\beta}R_{\alpha}$ for arbitrary 
   $\alpha,\beta\in L$.
\end{prop}
{\sc Proof}: The statement is trivial for the action on $\nu=0$.
So, let $\nu > 0$ be a (strictly) positive measure. We then have 
$\nu(X)=\|\nu\|$ and obtain
\begin{eqnarray*}
   \pi^{}_{\s\alpha}.\big((\pi^{}_{\s\alpha}.\nu)\otimes
     (\pi^{}_{\g\alpha}.\nu)\big) & = &
     \|\nu\|\,(\pi^{}_{\s\alpha}.\nu) \\
   \pi^{}_{\g\alpha}.\big((\pi^{}_{\s\alpha}.\nu)\otimes
     (\pi^{}_{\g\alpha}.\nu)\big) & = &
     \|\nu\|\,(\pi^{}_{\g\alpha}.\nu)
\end{eqnarray*}
in $\cM_+(X^{}_{\s\alpha})$ resp.\ $\cM_+(X^{}_{\g\alpha})$
where we adopt the same index convention for sets as we did for projectors.
Using $\| R^{}_{\alpha} (\nu)\|=\|\nu\|$ from Fact~\ref{positivity} and
the linearity of the mappings $\nu\mapsto (\pi.\nu)$, one can now apply 
the definition of the elementary recombinators to check explicitly that
$$ R^{}_{\alpha}\big(R^{}_{\alpha} (\nu)\big) \; = \; 
   R^{}_{\alpha}(\nu) \, .  $$

{}For commutativity, we may again assume $\nu >0$ and also
$\alpha < \beta$. Then 
\begin{eqnarray*}
  \pi^{}_{\s\alpha}.\big((\pi^{}_{\s\beta}.\nu)\otimes
    (\pi^{}_{\g\beta}.\nu)\big) & = & \|\nu\|\,(\pi^{}_{\s\alpha}.\nu) \\
  \pi^{}_{\g\alpha}.\big((\pi^{}_{\s\beta}.\nu)\otimes
    (\pi^{}_{\g\beta}.\nu)\big) & = &  \big(
    (\pi^{}_{\{\ceiling{\alpha},... ,\floor{\beta}\}}.\nu)\otimes
    (\pi^{}_{\g\beta}.\nu)\big) .
\end{eqnarray*}
The first equation can be verified directly, as in the previous case.
The second can easily be checked on Borel sets of the product form
$E=E^{}_{\{\ceiling{\alpha},... ,\floor{\beta}\}}\times
E^{}_{\g\beta}$, followed by an application of Fact~\ref{baire}.
Combining these intermediate results, one obtains 
$$ R^{}_{\alpha} \big( R^{}_{\beta} (\nu) \big) \; = \; 
   \frac{1}{\|\nu\|^2}\,\big((\pi^{}_{\s\alpha}.\nu) \otimes
   (\pi^{}_{\{\ceiling{\alpha},\dots ,\floor{\beta}\}}.\nu)\otimes
   (\pi^{}_{\g\beta}.\nu)\big) \; = \;
   R^{}_{\beta} \big( R^{}_{\alpha} (\nu) \big)\,,  $$
which proves our assertion.  \qed

\smallskip \noindent
{\sc Remark}: In view of positive homogeneity of the recombinators, see
Eq.~(\ref{pos-hom}), it would have been sufficient to prove our assertions
on $\cP(X)$. The above version, however, shows quite clearly where, and
how many, normalization factors $\|\nu\|$ appear in the tensor products. If we
restrict ourselves to probability measures below, one should keep this
in mind for extending arguments to the full cone, $\cM_+(X)$.

\smallskip
A close inspection of the proof of Proposition~\ref{comm-rules1} shows
that we have simultaneously proved the following useful property.
\begin{lemma} \label{forget}
   Let $\nu\in\cP(X)$ and $\alpha\in L$. For all
   $\beta\in L$ with $\beta\ge\alpha$, we have 
   $\pi^{}_{\s\alpha}. \big(R^{}_{\beta}(\nu)\big) =
   \pi^{}_{\s\alpha}.\nu$. Similarly, 
   $\pi^{}_{\g\alpha}. \big(R^{}_{\beta}(\nu)\big) =
   \pi^{}_{\g\alpha}.\nu$, for all $\beta\le\alpha$. \qed
\end{lemma}

\subsection{The IPL equation and its solution}

Let us start with a brief description of the recombination process for 
finite $X$, and a population of $m$ individuals, each of the form
$x = (x^{}_0, x^{}_1, \dots , x^{}_n )$ with $x_i \in X_i$.
Every individual carries a Poisson clock at each link $\alpha\in L$, with 
parameters $\varrho^{}_{\alpha} > 0$, which do not depend on the individual.
If the clock at link $\alpha$ of the individual $x$ rings, a random
partner $y$ is picked from the population for recombination at that link.
The recombined pair is then 
$(x^{}_0,\dots ,x^{}_{\floor{\alpha}},y^{}_{\ceiling{\alpha}},\dots, y^{}_n)$
and $(y^{}_0,\dots, y^{}_{\floor{\alpha}},x^{}_{\ceiling{\alpha}},\dots,x^{}_n)$.

To describe the entire population, let $Z_x(t)$ be the random variable that
gives the number of $x$-individuals at time $t$, and $Z(t)$ the combined
random vector with components $Z_x(t)$. Hence, if $Z(t)=z$, and $x\neq y$,
we can have transitions from $z$ to 
$z - u^{}_{x} - u^{}_{y} + u^{}_{(x^{}_{\s\alpha}\, ,\, y^{}_{\g\alpha})}
+ u^{}_{(y^{}_{\s\alpha}\, ,\, x^{}_{\g\alpha})}$, where we use our short 
hand notation for indices, and $u^{}_x$ to denote the unit vector 
corresponding to $x$. Such a transition occurs at rate 
$\varrho^{}_{\alpha} z^{}_{x} z^{}_{y}/(m- z^{}_x)$.

Note that this process
implies instant mixing of all (geno-)types in the population. This is an
idealization which neglects that maternal and paternal genes stay together
for the lifetime of an individual. Nevertheless, this is a good and realistic
model if recombination events are rare on the time scale of the
individual life span. This is certainly true if our sites
belong to the DNA sequence of a single gene, or a few adjacent genes. 
It is then well justified to describe recombination in terms of these
first order effects only.

Let us look at the influence of increasing $m$, whence we write
$Z^{(m)}(t)$ to indicate dependence on system size.
As $m \to \infty$, the sequence of random processes
$Z^{(m)}(t)/m$ converges almost surely to the solution of
a differential equation with initial condition
$Z^{(m)}(0)/m$ (resp.\ its limit as $m\to\infty$), 
see \cite[Thm.\ 11.2.1]{EK}. The corresponding IPL equation
\cite[Eq.~2.5]{EB}, reformulated in our measure-theoretic setting, 
reads 
\begin{equation} \label{r-master}
   \dot{\omega} \; = \; \Phi^{}_{\rm rec}(\omega) \; := \;
   \sum_{\alpha\in L} \varrho^{}_{\alpha} 
   \big(R_{\alpha}-\bs{1}\big) (\omega)\, .
\end{equation}
In line with our strategy for the mutation processes, we take this
nonlinear ODE as the general starting point for the recombination analysis
on product spaces $X$ built from arbitrary locally compact spaces $X_i$.
We will  assume that $\varrho^{}_{\alpha} > 0$, for all $\alpha\in L$,
without loss of generality (if $\varrho^{}_{\alpha} = 0$,
remove the link at $\alpha$, absorb the pair 
$(\floor{\alpha}, \ceiling{\alpha})$ into a single site, and identify 
$X_{\floor{\alpha}} \times X_{\ceiling{\alpha}}$ with the state
space at that site,
thus reducing the number of sites (and links) by one).

\begin{prop} \label{reco-dgl}
   The abstract Cauchy problem of the IPL equation $(\ref{r-master})$
   has a unique solution. Furthermore, $\cM_+(X)$ is positive
   invariant under the flow, with the norm of positive measures preserved. 
   In particular, $\cP(X)$ is positive invariant.
\end{prop}
{\sc Proof}: Consider $\dot{\omega}=\Phi^{}_{\rm rec}(\omega)$, which is a special
case of (\ref{gen-master}), so we want to apply Theorem~\ref{general-ODE}.   
By Fact \ref{lipschitz}, $\Phi^{}_{\rm rec}$ is Lipschitz, so assumption (A1) is
satisfied.

Let $\nu\in\cM_+(X)$, i.e.\ $\nu(E)\ge 0$ for all Borel sets $E\subset X$.
Let $E$ be any Borel subset of $X$ such that $\nu(E)=0$. Then
$$ \Phi^{}_{\rm rec}(\nu)(E) \; = \;
   \sum_{\alpha\in L} \varrho^{}_{\alpha} R_{\alpha}(\nu)(E) 
   \; \ge \;  0 $$
because each $R_{\alpha}(\nu)$ is a positive measure and all
$\varrho^{}_{\alpha} > 0$ by assumption, so (A2) is satisfied.

{}Finally, with Fact~\ref{positivity}, it is easy to check that 
$\Phi^{}_{\rm rec}(\nu)(X)=0$
for any positive measure $\nu$, which shows that assumption
(A3) is satisfied, too. Theorem \ref{general-ODE} then
establishes our claims. \qed

\smallskip
The difficulty in solving $(\ref{r-master})$ stems from the nonlinearity
of the right-hand side, so $\Phi^{}_{\rm rec}$ cannot be considered as the 
generator
of an exponential semigroup. It is, however, rather natural to expect
that the solution should still have a rather similar structure, as the
$R_{\alpha}$ are at least positive homogeneous of degree one and
commute with one another. Let us therefore, for any $G\subset L$, introduce
the {\em composite} recombinators
\begin{equation} \label{composite}
   R^{}_G \; := \; \prod_{\alpha\in G}R_{\alpha}\, .
\end{equation}
They are well-defined on $\cM_+(X)$ due to Proposition \ref{comm-rules1}, 
while an order of the product has to be specified otherwise. In any case,
$\|R^{}_G (\omega)\|\le\|\omega\|$ for all $\omega\in\cM(X)$.
Note that $R^{}_{\es}=1$ and $R^{}_{\{\alpha\}} = R^{}_{\alpha}$
in this notation. The composite recombinators are again positive
homogeneous of degree one. A simple induction argument based on
Proposition~\ref{comm-rules1} gives the following result.
\begin{coro} \label{alg-rules}
   On $\cM_+(X)$, the composite recombinators satisfy 
   $$ R^{}_{G} R^{}_{H} \; = \; R^{}_{G\cup H} \, , $$ 
   for arbitrary $G,H\subset L$.  Furthermore, each $R^{}_G$
   maps $\cM_+(X)$ into itself and preserves the norm of positive
   measures. \qed
\end{coro}

Let us pretend for a moment that the idempotents $R^{}_{\alpha}$ were 
actually linear operators. In such a case, we would get
$$ \exp\bigl( \varrho^{}_{\alpha} t (R^{}_{\alpha}-\bs{1})\bigr)
  \; = \; \exp(-\varrho^{}_{\alpha} t ) \bs{1} +
  \bigl(1 - \exp(-\varrho^{}_{\alpha} t)\bigr) R^{}_{\alpha}\,. $$
Taking the product over such terms for all $\alpha\in L$ and expanding
it would formally lead to the sum
$$ \sum_{G\subset L} a^{}_{G}(t)\, R^{}_G $$
with the coefficient functions
\begin{equation} \label{coeff-fun}
  a^{}_{G}(t) \; = \; \Big(
      \prod_{\alpha\in\overline{G}}\exp(-\varrho^{}_{\alpha} t)\Big)
      \cdot \Big(
      \prod_{\beta\in G}\, \bigl(1 - \exp(-\varrho^{}_{\beta} t)\bigr) 
      \Big).
\end{equation}
It will have a touch of magic below when we prove that this little
``derivation'' actually gives the correct answer! After we have established
our main result in Theorem~\ref{sol1}, we will come back to these coefficients
and give them a probabilistic interpretation. This will also motivate
why they are a very reasonable guess to start with.

As mentioned before, the elementary recombinators are {\em not\/}
linear. Nevertheless, they have a related property on convex combinations.
If $\omega=\sum_{i=1}^k a^{}_i \nu^{}_i$ is a convex linear combination
of positive measures $\nu^{}_i$ of equal norm, we get
\begin{equation} \label{non-linear}
   R^{}_{\alpha}(\omega) \; = \; 
     \sum_{i=1}^k a^{}_i \, R^{}_{\alpha}(\nu^{}_i)
     \, + \, {\mathfrak R}
\end{equation}
where one can show, by a rather straight-forward calculation which we omit
here, that the remainder $\mathfrak R$ is given by
$${\mathfrak R} \; = \; - \, \sum_{i\s j} a^{}_i a^{}_j \, 
   \frac{\|\nu^{}_i - \nu^{}_j\|}{\|\omega\|}
   \, R^{}_{\alpha}(\nu^{}_i - \nu^{}_j)\, .$$
This shows that the recombinators are indeed inherently
nonlinear, but also that they might act like linear operators on 
{\em special\/} convex combinations, namely those for which the 
remainder vanishes. This is precisely what we need to solve our problem.
\begin{prop}  \label{conv-linear}
   Let\/ $\nu$ be a positive measure, $\alpha\in L$ and 
   $a^{}_G(t)$ the coefficient functions of $(\ref{coeff-fun})$. Then,
   for any fixed $t\ge 0$, we have 
   $$R^{}_{\alpha} \Big( \sum_{G\subset L} a^{}_G(t)\, R^{}_G (\nu) \Big) \; = \; 
   \sum_{G\subset L} a^{}_{G}(t)\, R^{}_{G\cup\{\alpha\}} (\nu)\, .$$
\end{prop}

Before we prove this result, we formulate a special property of the coefficient
functions first. Observe that, for fixed $t\ge 0$, 
$q^{}_{\alpha} := \exp(-\varrho^{}_{\alpha} t)$ is a number between 
$0$ and $1$. It can be interpreted as a probability (namely that
link $\alpha$ has not been hit until time $t$). With this, the
coefficients read
\begin{equation} \label{temp-coeff}
    a^L_G \; = \; \prod_{\alpha\in L-G} q^{}_{\alpha}\;
    \prod_{\beta\in G} (1-q^{}_{\beta}) \; \ge \; 0
\end{equation}
where we have suppressed the (fixed) time, but added the set of links, $L$, 
as an upper index. We can now formulate a crucial factorization property.
\begin{lemma} \label{coeff-factorization}
  Let $L=L_1 \dot{\cup} L_2$ be a partition of $L$, and set $G_i = G\cap L_i$
  for an arbitrary $G\subset L$. Then, the coefficients of $(\ref{temp-coeff})$
  satisfy $a^L_G = a^{L_1}_{G_1} \cdot a^{L_2}_{G_2}$.
  Furthermore, for any $L^{'}\subset L$, we have 
  $$\sum_{H\subset L^{'}} a^{L^{'}}_H  \; = \; 1\, .$$
\end{lemma}
{\sc Proof}: Since $L_1\cap L_2 = \varnothing$, the first statement is a
direct consequence of the product form of $a^L_G$ in Eq.~(\ref{temp-coeff}).
The normalization property can be verified from the probabilistic
interpretation mentioned above. If $1-q^{}_{\alpha}$ (resp.\ $q^{}_{\alpha}$)
is the probability that link $\alpha$ has (resp.\ has not) been hit,
$a^{L^{'}}_H$ is the probability that, of the links in $L^{'}$,
precisely $H$ is spared. Consequently,
$\sum_{H\subset L^{'}} a^{L^{'}}_H$ is the sum over the probabilities of all
possible events, hence equal to 1. Alternatively, this identity
can be derived from a simple M\"obius inversion argument, as we show
below in Fact~\ref{conv-comb}. \qed
 
\smallskip\noindent
{\sc Proof of Proposition~\ref{conv-linear}}: 
Since the recombinators are positive homogeneous of degree one,
it suffices to prove the statement for $\nu$ a probability measure.
Let $\alpha\in L$ be fixed.

Set $\omega = \sum_{G\subset L} a^L_G \, R^{}_G (\nu)$. 
Since $\nu\in {\mathcal P}(X)$ implies 
$\omega \in {\mathcal P}(X)$, we obtain
\begin{eqnarray*}
   R^{}_{\alpha}(\omega)  & = &  R^{}_{\alpha}
   \Big( \sum_{G\subset L} a^L_G \, R^{}_G (\nu) \Big) \\
   & = & \bigg( \pi^{}_{\s\alpha} . \Big( \sum_{G\subset L}
   a^L_G \, R^{}_G (\nu) \Big) \bigg) \otimes
   \bigg( \pi^{}_{\g\alpha} . \Big( \sum_{H\subset L}
   a^L_H \, R^{}_H (\nu) \Big) \bigg)   \\
   & = & \sum_{G,H \subset L} a^L_G \, a^L_H \,
   \Big( \big(\pi^{}_{\s\alpha} . R^{}_G (\nu)\big) \otimes
   \big(\pi^{}_{\g\alpha} . R^{}_H (\nu)\big)\Big)
\end{eqnarray*}
where we have used the linearity of the mappings $\pi^{}_{\s\alpha}$
and $\pi^{}_{\g\alpha}$.

Let us define $L_1 = \{\frac{1}{2}, \frac{3}{2}, \ldots, \alpha \}$
and $L_2 = L - L_1$, so that $L = L_1 \dot{\cup} L_2$ is a partition
of $L$. Also, let $G_i = G \cap L_i$ and $H_i = H \cap L_i$, for
$G,H \subset L$. Lemma~\ref{forget} then tells us that
\begin{eqnarray*}
  \big(\pi^{}_{\s\alpha} . R^{}_G (\nu)\big) \otimes
   \big(\pi^{}_{\g\alpha} . R^{}_H (\nu)\big)
  & = & \big(\pi^{}_{\s\alpha} . R^{}_{G_1} (\nu)\big) \otimes
   \big(\pi^{}_{\g\alpha} . R^{}_{H_2} (\nu)\big) \\
  & = &  \big(\pi^{}_{\s\alpha} . R^{}_{G_1 \cup H_2} (\nu)\big) \otimes
   \big(\pi^{}_{\g\alpha} . R^{}_{G_1 \cup H_2} (\nu)\big) \\
  & = & R^{}_{\alpha} \big( R^{}_{G_1 \cup H_2} (\nu) \big).
\end{eqnarray*}
Inserting this into the previous equation and invoking 
Lemma~\ref{coeff-factorization} repeatedly gives
\begin{eqnarray*}
   R^{}_{\alpha} (\omega) & = & \sum_{G,H \subset L}
   a^L_G  \, a^L_H  \, 
   R^{}_{\alpha} \big( R^{}_{G_1 \cup H_2} (\nu) \big)  \\
   & = & \sum_{G_1 \subset L_1} \sum_{G_2 \subset L_2}
   \sum_{H_1 \subset L_1} \sum_{H_2 \subset L_2} 
   a^{L_1}_{G_1} \, a^{L_2}_{G_2} \, 
   a^{L_1}_{H_1} \, a^{L_2}_{H_2} \,
   R^{}_{\alpha} \big( R^{}_{G_1 \cup H_2} (\nu) \big)  \\
   & = & \sum_{G_1 \subset L_1} \sum_{H_2 \subset L_2}
   a^{L_1}_{G_1} \, a^{L_2}_{H_2} \,
   R^{}_{\alpha} \big( R^{}_{G_1 \cup H_2} (\nu) \big)  \\
   & = & \sum_{K\subset L} a^L_K \,
   R^{}_{\alpha} \big( R^{}_{K} (\nu) \big)\,, 
\end{eqnarray*}
which proves our assertion.   \qed

\smallskip \noindent
{\sc Remark}: Proposition \ref{conv-linear} admits the following interpretation.
Let $\nu$ be a positive measure, with $\|\nu\|=m>0$. Then, the $2^{|L|}$ 
measures $R^{}_G (\nu)$ with $G\subset L$ form the vertices of a 
$\|.\|$-closed simplex in $\cM_+^m(X)$. On some of their convex combinations
(in particular along solutions, as we will see shortly), the elementary recombinators 
$R^{}_{\alpha}$ act linearly. It is this simplex, foliated into solution
curves, to which the entire time evolution is constrained, with
$\nu$ as the initial condition.

The positive measure in Proposition \ref{conv-linear} was arbitrary.
This means that, when restricting the action of the $R^{}_{\alpha}$'s
to $\cM_+(X)$, we can formulate the rule on the level of operators. 
Observe that 
$R_{\alpha} \big(R_K (\nu)\big) = R_{K\cup\{\alpha\}}(\nu) = R_K (\nu')$ 
where $\nu' = R_{\alpha}(\nu)$.
By a simple induction argument, we thus arrive at
\begin{coro} \label{new-rules}
   Let\/ $a^{}_G (t)$ be the coefficient function of $(\ref{coeff-fun})$,
   and let $t\ge 0$ be fixed. 
   On $\cM_+(X)$, the recombinators satisfy the equation
$$ R^{}_H \Big( \sum_{G\subset L} a^{}_G (t) \, R^{}_G \Big) \; = \;
   \sum_{G\subset L} a^{}_G (t) \, R^{}_{G\cup H} $$
   for arbitrary $H \subset L$.  \qed
\end{coro}

We now assume that the initial condition, $\omega^{}_0$, is a positive measure
and make the following {\em ansatz\/} for the solution of (\ref{r-master}): 
\begin{equation} \label{ansatz}
   \omega^{}_t \; = \; \sum_{G\subset L} a^{}_{G}(t)\, R^{}_G(\omega^{}_0) 
\end{equation}
with the coefficient functions $a^{}_{G}(t)$ of (\ref{coeff-fun}). Note
that they do {\em not\/} depend on $\omega^{}_0$.
The initial values are $a^{}_{\es}(0)=1$ and 
$a^{}_{G}(0)=0$ for all $\es\neq G\subset L$.
By Corollary \ref{alg-rules}, each $R^{}_G(\omega^{}_0)$ is a positive measure 
with the same norm as $\omega^{}_0$. This implies that, as long as
$a^{}_{G}(t)\ge 0$ for $t\ge 0$, the ansatz for $\omega^{}_t$ must form a
convex linear combination of positive measures of equal norm if it
is a solution of (\ref{m-master}). This follows from Eq.~(\ref{temp-coeff})
together with Lemma~\ref{coeff-factorization}, or from Fact~\ref{conv-comb} below.

The time derivative of $\omega^{}_t$ of (\ref{ansatz}) is $\dot{\omega}^{}_t = 
\sum_{G\subset L} \dot{a}^{}_{G}(t)\, R^{}_G(\omega^{}_0)$. On
the other hand, Proposition \ref{conv-linear} means that the 
$R^{}_{\alpha}$ act linearly on the convex combination (\ref{ansatz}), 
and we obtain 
\begin{eqnarray*}
  \Phi^{}_{\rm rec} (\omega^{}_t) & = & \sum_{\alpha\in L} \varrho^{}_{\alpha}
  \, (R^{}_{\alpha} - \bs{1}) (\omega^{}_t)  \\
  & = & \sum_{\alpha\in L} \varrho^{}_{\alpha} 
        \sum_{G\subset L} a^{}_{G} (t)\,
        \big(R^{}_{G\cup\{\alpha\}}(\omega^{}_0) 
         - R^{}_{G}(\omega^{}_0)\big) \\
  & = & \sum_{\alpha\in L} \varrho^{}_{\alpha} \Big[
        \sum_{\alpha\in \text{\d{$G$}}\subset L} a^{}_{G\setminus\{\alpha\}} (t)
              \, R^{}_{G}(\omega^{}_0) \; -  
        \sum_{\alpha\not\in \text{\d{$G$}}\subset L} a^{}_{G} (t)
              \, R^{}_{G}(\omega^{}_0) \, \Big] \\
  & = & \sum_{G\subset L} \Big[ \sum_{\alpha\in G} \varrho^{}_{\alpha}\,
        a^{}_{G\setminus\{\alpha\}}(t) \, - \, \sum_{\beta\in\overline{G}}
        \varrho^{}_{\beta}\, a^{}_{G}(t) \Big] R^{}_{G} (\omega^{}_0)\, ,
\end{eqnarray*}
where we use the notation {\d{$G$}} in the third step to indicate the summation 
variable. It is now a straight-forward calculation to check that the coefficients
$a^{}_{G}(t)$ of (\ref{coeff-fun}) indeed satisfy the equations
$$ \dot{a}^{}_{G}(t) \; = \; \sum_{\alpha\in G} \varrho^{}_{\alpha}\,
        a^{}_{G\setminus\{\alpha\}}(t) \, - \, \sum_{\beta\in\overline{G}}
        \varrho^{}_{\beta}\, a^{}_{G}(t) $$
and that they constitute a convex combination in (\ref{ansatz}).
Consequently, our ansatz solves the IPL equation (\ref{r-master}),
and, by Proposition \ref{reco-dgl}, this is the unique solution we are after. 
We have thus established the following main result.
\begin{theorem} \label{sol1}
   The ansatz $(\ref{ansatz})$ solves the IPL equation $(\ref{r-master})$ 
   with initial condition\/ $\omega^{}_0 \in\cM_+(X)$ if and only if the 
   coefficient functions are given by $\,(\ref{coeff-fun})$, i.e.\ by
   $$ a^{}_{G}(t) \; = \; 
      \exp\Bigl(-\sum_{\alpha\in\overline{G}}\,\varrho^{}_{\alpha} t\Bigr)
      \cdot \prod_{\beta\in G}\, \bigl(1 - \exp(-\varrho^{}_{\beta} t)\bigr) $$
   for all\/ $G\subset L$. \qed
\end{theorem}
\noindent {\sc Remark}:
To interpret the coefficient $a^{}_G (t)$, let us consider a
single individual.
Since $\exp(-\varrho^{}_{\alpha} t)$ is the probability
that  link $\alpha$ has experienced no crossover event until
time $t$ (recall that we have assumed a Poisson process of
rate $\varrho^{}_{\alpha}$ at link $\alpha$), $a^{}_{G}(t)$ 
may be interpreted as  the probability
that the set of all links that have, up to time $t$, experienced
at least one crossover event, is precisely $G$.

Note that the above result relies on the assumption of single,
independent crossover events, which is described by recombinators
that commute. In more general models, with  multiple, dependent
events, the algebraic structure is rather involved. This makes
solutions much more cumbersome, or considerably less
explicit in nature (for review, see \cite[Ch.~6]{Lyu}).  

Let us come back to the meaning of Eq.~(\ref{ansatz}) in combination with
Theorem~\ref{sol1}. If $\varphi^{}_t$ denotes the flow of the IPL equation
(\ref{r-master}), we obtain, for all $t\ge 0$, the identity
\begin{equation} \label{flow-def}
  \varphi^{}_t \; = \; \sum_{G\subset L} a^{}_G(t)\, R^{}_G 
\end{equation}
which is valid on the cone $\cM_+(X)$. As usual, $\varphi^{}_0=\bs{1}$
and $\varphi^{}_t \circ \varphi^{}_s = \varphi^{}_{t+s}$, for all
$t,s \ge 0$. This implies the identity
$$ a^{}_G (t+s) \; = \; \sum_{\substack{H,K\subset L\\H\cup K = G}}
   a^{}_H (t)\, a^{}_K (s)\,,  $$
which can be verified by direct computation. More interestingly, we
also have
\begin{fact} \label{flow}
   On $\cM_+(X)$, the forward flow of\/ $(\ref{r-master})$ commutes with the
   recombinators, i.e.\ $R^{}_G \circ \varphi^{}_t = \varphi^{}_t \circ R^{}_G$,
   for all $t\ge 0$ and $G\subset L$.
\end{fact}
{\sc Proof}: Let $\nu\in\cM_+(X)$ and fix $G\subset L$. Then
\begin{eqnarray*}
   R^{}_G \big(\varphi^{}_t (\nu)\big) & = &
   R^{}_G \Big( \sum_{H\subset L} a^{}_H (t) R^{}_H (\nu) \Big)
   \,\; = \,\; \sum_{H\subset L} a^{}_H (t) R^{}_{G\cup H} (\nu) \\
   & = & \sum_{H\subset L} a^{}_H (t) R^{}_H \big( R^{}_G (\nu)\big)
   \,\; = \,\; \varphi^{}_t \big( R^{}_G (\nu)\big)
\end{eqnarray*}
by an application of Corollary \ref{new-rules}.  \qed

\smallskip
Once the solution is known, the remaining task is to identify
linear combinations of the $R^{}_H(\omega)$ that decouple from each
other and decay exponentially.
To this end, we employ combinatorial techniques to regroup
the terms of the solution according to their exponential
damping factors. Let us first expand the expression for $a^{}_{G}(t)$,
$$ a^{}_{G}(t) \; = \; \sum_{K\subset G} (-1)^{|G-K|}\,
   \exp\Bigl(-\sum_{\alpha\in\overline{K}}\,\varrho^{}_{\alpha} t\Bigr)\, .$$
This suggests to define new functions $b^{}_K(t)$ via 
\begin{equation} \label{new-coeff}
   b^{}_K(t) \; = \; 
 \exp\bigl(-\sum_{\alpha\in\overline{K}}\,\varrho^{}_{\alpha} t\bigr) \, ,
\end{equation}
with the usual convention that the empty sum is 0.
In particular, we have $b^{}_{\es}(t)=a^{}_{\es}(t)=
\exp\bigl(-\sum_{\alpha\in L}\,\varrho^{}_{\alpha} t\bigr)$ and
$b^{}_{K}(0) = 1$ for all $K\subset L$.
Now, the M\"obius inversion of (\ref{moebius1}) and (\ref{moebius2}), 
used  backwards, gives us the relation
$$ b^{}_K(t) \; = \; \sum_{G\subset K} a^{}_G(t) \, .$$
One immediate consequence is
\begin{equation} \label{conv-sum-check}
  \sum_{G\subset L} a^{}_{G}(t) \; = \; b^{}_{L}(t) \; = \; 1 \, .
\end{equation}
So, together with the observation that the functions $a^{}_{G}(t)$ of
Theorem \ref{sol1} are always non-negative, we have independently confirmed
\begin{fact} \label{conv-comb}
   If\/ $\omega^{}_0\in\cM_+(X)$, the coefficient functions\/ 
   $a^{}_{G}(t)$ of Theorem $\ref{sol1}$ constitute a convex linear 
   combination of positive measures in Eq.~$(\ref{ansatz})$. \qed
\end{fact}

The significance of the new functions becomes clear by realizing that
there is an analogue on the level of operators. To this end, we rewrite
the composite recombinators in terms of new operators via 
$R^{}_H = \sum_{\text{\d{$G$}}\supset H} T^{}_G$ and obtain, by an 
obvious variant of  M\"obius inversion,
\begin{equation} \label{t-def}
   T^{}_G \; := \; \sum_{\text{\d{$H$}}\supset G} (-1)^{|H-G|}\, R^{}_H \, .
\end{equation}
A straight-forward calculation then reveals that
\begin{equation} \label{new-form}
  \omega^{}_t \; = \; \sum_{G\subset L} a^{}_{G}(t)\, R^{}_G(\omega^{}_0) 
       \; = \; \sum_{K\subset L} b^{}_{K}(t)\, T^{}_K(\omega^{}_0)\, . 
\end{equation}
Note that, as a consequence of Eqs.~(\ref{pos-hom}), (\ref{composite})
and (\ref{t-def}), the operators $T^{}_G$ are positive homogeneous of
degree one, i.e.
\begin{equation} \label{t-hom}
    T^{}_G(\alpha \nu) \; = \; | \alpha | \cdot T^{}_G(\nu)\,.
\end{equation}

Let us now introduce new measures 
$\nu^{}_{G}(t) := b^{}_{G}(t) \, T^{}_{G} (\omega^{}_0)$, which are
elements of $\cM(X)$, but no longer positive in general.
\begin{prop} \label{diseq1}
   The signed measures\/ $\nu^{}_{G}(t)$ solve the Cauchy problem
$$ \dot{\nu}^{}_G (t) \; = \;
   - \Bigl(\mbox{$\sum_{\alpha\in \overline{G}}$ } 
    \varrho^{}_{\alpha}\Bigr) \, \nu^{}_{G}(t) $$
   with initial condition $\nu^{}_{G}(0) = T^{}_{G} (\omega^{}_0)$,
   for all $G\subset L$.
\end{prop}
{\sc Proof}: The result is a direct consequence of the fact that the
coefficient functions $b^{}_{G}(t)$ solve the ordinary initial value problems
$$\dot{b}^{}_G (t) \; = \;
   - \Bigl(\mbox{$\sum_{\alpha\in \overline{G}}$} \, \varrho^{}_{\alpha}\Bigr)
   \, b^{}_{G}(t) $$ 
with initial conditions $b^{}_{G}(0) = 1$, see above. \qed

\smallskip
So, the transformation (\ref{t-def}) resulted in regrouping the terms
of the solution to the IPL equation (\ref{r-master}) according to their 
exponential decay factors in time. In particular,
$$ \nu^{}_{L} (t) \; \equiv \; T^{}_{L}(\omega^{}_0)
   \; = \; R^{}_{L}(\omega^{}_0) \; = \; 
   \bigotimes_{i=0}^{n}\, (\pi^{}_i . \omega^{}_0) $$
is the unique {\em limit measure\/} of the process starting from $\omega^{}_0$. 
Due to the action of $R^{}_{L}$, it is a complete product measure and 
reflects total independence, and we obtain $\omega^{}_t \to R^{}_{L}(\omega^{}_0)$
as $t\to\infty$ in the $\|.\|$-topology. This is so because
$$ \|\omega^{}_{t} - \nu^{}_{L}(t)\| \; = \;
   \Big\|\sum_{K\subsetneq L} b^{}_{K} (t)\, T^{}_{K}(\omega^{}_{0})\Big\| 
   \; \le \; \sum_{K\subsetneq L} b^{}_{K} (t)\, \|T^{}_{K}(\omega^{}_{0})\| $$
where all remaining coefficient functions $b^{}_{K} (t)$, i.e.\ those with
$K\subsetneq L$, decay exponentially
(recall that $\varrho^{}_{\alpha} > 0$ for all $\alpha\in L$).

\subsection{Linkage disequilibria}

Starting from the measures $\nu^{}_{G}(t)$, we will now  identify a minimal, 
complete set of variables by evaluating
certain $k$-point cylinder functions (called $k$-point functions from now
on) or correlation functions known as linkage disequilibria
in genetics. 
They are important for data analysis because they allow 
to evaluate associations between sites up to a given order
from measured type frequencies, and average over all others by
marginalization.
This way, a certain amount of stochasticity, which is present
in all real (finite) populations, is smoothed out.

Various different definitions of linkage disequilibria
are available in the literature (see \cite[p.\ 183--186]{Reinhard} for 
an overview). But only special choices decouple (see \cite{Dawson,Dawson2}),
and these are the linkage disequilibria we are after.
In view of the applications, 
we will now restrict ourselves to the case that $X$ is a {\em finite\/} set, 
although the results hold, with only minor modifications, also more generally. 
Eq.\ (\ref{new-form}) and Proposition \ref{diseq1} suggest to employ
the signed measures $T^{}_{G}(\omega^{}_{0})$. The corresponding functions
$b^{}_{G}(t)$ will then describe their evolution in time.

Let $\langle j^{}_1,\dots ,j^{}_k \rangle$, with $j^{}_1 < \dots < j^{}_k$, 
symbolically denote a cylinder set 
in $X=X^{}_N$ which is specified at sites $j^{}_i$, for $1\le i \le k$. More
specifically, these are sets of the product form
$$\langle j^{}_1,\dots ,j^{}_k \rangle \; = \;
  X^{}_{\{0,...,j^{}_1 - 1\}} \times 
   \{x^{}_{j^{}_1}\} \times \big[...\big] \times  \{x^{}_{j^{}_k}\} 
  \times X^{}_{\{j^{}_k + 1, ... , n\}} $$
where $[...]$ contains factors $\{x^{}_i\}$ or $X_i$ depending on
whether $i$ appears in $\langle j^{}_1, ... \, ,j^{}_k \rangle$ or not.
For $\nu\in\cP(X)$ and arbitrary $\alpha\in L$, we then have
\begin{eqnarray*}
   R^{}_{\alpha} (\nu) \big(\langle j^{}_1,...\, ,j^{}_k \rangle\big)
   & = & \big((\pi^{}_{\s\alpha}.\nu)\otimes (\pi^{}_{\g\alpha}.\nu)\big)
         \big(\langle j^{}_1,...\, ,j^{}_k \rangle\big)  \\[2mm]
   & = & \begin{cases}
          \nu\big(\langle j^{}_1,...\, ,j^{}_k \rangle\big) &
           \text{if $\alpha < j^{}_1$ or $\alpha > j^{}_k$} \\
          \nu\big(\langle j^{}_1,...\, ,j^{}_s \rangle\big) \,
          \nu\big(\langle j^{}_{s+1},...\, ,j^{}_k \rangle\big) &
           \text{if $j^{}_s < \alpha < j^{}_{s+1}\,$.}
         \end{cases}
\end{eqnarray*}
For later convenience, we also define $\langle\es\rangle = X$
so that $R^{}_{\alpha}(\nu)\big(\langle\es\rangle\big) = 1$.
\begin{lemma} \label{zero-rule}
   If\/ $\nu\in\cP(X)$, we have\/ $T^{}_{G} (\nu) 
   \big(\langle j^{}_1,\dots ,j^{}_k \rangle\big) = 0$
   whenever the set\/ $\overline{G}$ contains an element 
   that is less than $j^{}_1$ or larger than $j^{}_k$.
\end{lemma}
{\sc Proof}: Let $I=\{\beta\mid j^{}_1 < \beta < j^{}_k\}$.
Assume there is an $\alpha\in\overline{G\cup I}$. Then 
\begin{eqnarray*}
  \lefteqn{T^{}_{G} (\nu) \big(\langle j^{}_1,\dots ,j^{}_k 
           \rangle\big) } \\[2mm]
   & = & \sum_{\text{\d{$H$}}\supset G} (-1)^{|H-G|} R^{}_{H}(\nu)
         \big(\langle j^{}_1,\dots ,j^{}_k \rangle\big) \\
   & = & \sum_{\text{\d{$H$}}\supset G} (-1)^{|H-G|} R^{}_{H\cap I}(\nu) 
         \big(\langle j^{}_1,\dots ,j^{}_k \rangle\big) \\
   & = & \sum_{\alpha\not\in\text{\d{$H$}}\supset G}
         \big( (-1)^{|H-G|} + (-1)^{|(H\cup\{\alpha\})-G|} \big)
         R^{}_{H\cap I}(\nu) 
         \big(\langle j^{}_1,\dots ,j^{}_k \rangle\big) 
\end{eqnarray*}
where the previous calculation was used in the second step, and
summation is over $H$. Clearly, the last expression vanishes. \qed

\smallskip
Let us now define the time-dependent $k$-point functions  as
\begin{equation} \label{corr}
   F_G^{(t)} (j^{}_1,\dots ,j^{}_k) \; = \;
    T^{}_G (\omega^{}_t) 
   \big(\langle j^{}_1,\dots ,j^{}_k \rangle \big)
\end{equation} 
for arbitrary $G\subset L$, where the notation is again symbolic in that
we only specify the positions $j^{}_i$, but not the corresponding values.
To relate this to Eq.~(\ref{new-form}), we show 
\begin{prop} \label{measurement}
   If\/ $\omega^{}_0 \in\cM_+(X)$, we have
   $ T^{}_G (\omega^{}_t)= b^{}_G (t)\, T^{}_G (\omega^{}_0)$,
   for all $G\subset L$ and $t\ge 0$.
\end{prop}
{\sc Proof}:
Since $b^{}_G(0)=1$, equality holds for $t=0$, and the claim follows if
we show that  $T^{}_G (\omega^{}_t)$ and $b^{}_G(t)\, T^{}_G (\omega^{}_0)$
satisfy the same differential equation. With $\omega^{}_t =
\varphi^{}_t (\omega^{}_0)$, compare (\ref{flow-def}), we obtain
\begin{eqnarray*}
  \frac{d}{dt} T^{}_G (\omega^{}_t) & = & \frac{d}{dt}
  \sum_{\text{\d{$H$}}\supset G} (-1)^{|H-G|} R^{}_H 
  \big( \varphi^{}_t (\omega^{}_0) \big)  \\
  & = & \sum_{\text{\d{$H$}}\supset G} (-1)^{|H-G|}\, \frac{d}{dt}
  \varphi^{}_t \big( R^{}_H (\omega^{}_0)\big) 
  \qquad \qquad \qquad \quad \mbox{(by Fact \ref{flow})} \\
  & = & \sum_{\text{\d{$H$}}\supset G} (-1)^{|H-G|}\,
  \Phi^{}_{\rm rec} \big(\varphi^{}_t ( R^{}_H (\omega^{}_0))\big)  
  \qquad \qquad \quad \; \mbox{(by Eq.~(\ref{r-master}))} \\
  & = & \sum_{\text{\d{$H$}}\supset G} (-1)^{|H-G|}
  \sum_{\alpha\in \overline{G}} \varrho^{}_{\alpha} \big(R^{}_{H\cup \{\alpha\}}
  - R^{}_H\big) (\omega^{}_t) 
   \qquad \; \mbox{(by Fact \ref{flow})} \\
  & = & - \, \Big(\sum_{\alpha\in\overline{G}}\varrho^{}_{\alpha}\Big) 
   \, T^{}_G (\omega^{}_t) \, .
\end{eqnarray*}
The last step is correct because
$$ \sum_{\alpha\in \overline{G}} \varrho^{}_{\alpha}
   \sum_{\text{\d{$H$}}\supset G} (-1)^{|H-G|} R^{}_{H\cup \{\alpha\}} 
   \; = \; 0 $$
by an argument analogous to the one used in the last
step of the proof of Lemma~\ref{zero-rule}. Now, a comparison with
Proposition~\ref{diseq1} establishes the claim.   \qed

\smallskip
Even after Lemma \ref{zero-rule}, there are still too many functions 
around. It is thus reasonable to select an independent set from them.
To see how to do this, assume that we have an index $\alpha\in G\cap I$,
with $I=\{\beta\mid j^{}_1 < \beta < j^{}_k\}$ for a cylinder set of
type $\langle j^{}_1,\dots ,j^{}_k \rangle$ as above.
Let $H$ be a subset of $L$ that contains $G$, so $\alpha\in H$ in
particular, and $\nu \in \cP(X)$. Then $R^{}_H (\nu) = 
R^{}_{\alpha} (\nu^{}_H)$ with
$\nu^{}_H = R^{}_{H\setminus\{\alpha\}} (\nu)$.
The little calculation before Lemma \ref{zero-rule} now tells us that
\begin{eqnarray*}
  R^{}_{H} (\nu) \big(\langle j^{}_1,\dots ,j^{}_k \rangle\big) 
  & = & R^{}_{\alpha} (\nu^{}_H)
        \big(\langle j^{}_1,\dots ,j^{}_k \rangle\big)  \\
  & = & \big[ \nu^{}_H \big(\langle j^{}_1,\dots ,j^{}_s \rangle\big) 
        \big]\!\cdot\!
        \big[ \nu^{}_H \big(\langle j^{}_{s+1},\dots ,j^{}_k \rangle
        \big) \big]  \\
  & = & \big[ R^{}_{\alpha}(\nu^{}_H) 
        \big(\langle j^{}_1,\dots ,j^{}_s \rangle\big) \big]\!\cdot\!
        \big[ R^{}_{\alpha}(\nu^{}_H) 
        \big(\langle j^{}_{s+1},\dots ,j^{}_k \rangle \big) \big]  \\
  & = & \big[ R^{}_{H} (\nu)
        \big(\langle j^{}_1,\dots ,j^{}_s \rangle\big) \big]\!\cdot\!
        \big[ R^{}_{H} (\nu)  
        \big(\langle j^{}_{s+1},\dots ,j^{}_k \rangle \big) \big]
\end{eqnarray*}
where $j^{}_s < \alpha < j^{}_{s+1}$. Consequently, defining
$I_1 = \{\beta\mid j^{}_1 < \beta < j^{}_s \}$ and
$I_2 = \{\beta\mid j^{}_{s+1} < \beta < j^{}_k \}$, and 
referring back to (\ref{t-def}), we also get
\begin{eqnarray*}
    \lefteqn{T^{}_{G} (\nu)\big(\langle 
             j^{}_1,\dots ,j^{}_k \rangle\big) }  \\ [2mm]
    & = & \sum_{\text{\d{$H$}}\supset G} (-1)^{|H-G|}\;
          \big[ R^{}_{H} (\nu)
          \big(\langle j^{}_1,\dots ,j^{}_s \rangle\big) \big]\!\cdot\!
          \big[ R^{}_{H} (\nu)  
          \big(\langle j^{}_{s+1},\dots ,j^{}_k \rangle \big) \big] \\
   & = &  \sum_{\text{\d{$H$}}\supset G} (-1)^{|H-G|} 
          \sum_{\substack{\text{\d{$K\!$}}_1 \supset H \\
                          \text{\d{$K\!$}}_2 \supset H } }
          \big[ T^{}_{K_1} (\nu)
          \big(\langle j^{}_1,\dots ,j^{}_s \rangle\big) \big]\!\cdot\!
          \big[ T^{}_{K_2} (\nu)  
          \big(\langle j^{}_{s+1},\dots ,j^{}_k \rangle \big) \big] \\
   & = &  \sum_{\text{\d{$H$}}\supset G} (-1)^{|H-G|} \!\!\!\!\!
          \sum_{\substack{\text{\d{$K\!$}}_1 \supset (H\cup\overline{I}_1) \\
                          \text{\d{$K\!$}}_2 \supset (H\cup\overline{I}_2) } }
          \!\!\!\!\! \big[ T^{}_{K_1} (\nu)
          \big(\langle j^{}_1,\dots ,j^{}_s \rangle\big) \big]\!\cdot\!
          \big[ T^{}_{K_2} (\nu)  
          \big(\langle j^{}_{s+1},\dots ,j^{}_k \rangle \big) \big]
\end{eqnarray*} 
where Lemma~\ref{zero-rule} was used in the last step to remove
terms that vanish.
This equation means that $T^{}_{G} (\nu)\big(\langle 
j^{}_1,\dots ,j^{}_k \rangle\big)$, whenever an $\alpha\in G\cap I$ 
exists, either vanishes (if Lemma~\ref{zero-rule} applies) or is
a polynomial expression in $\ell$-point functions with $\ell < k$.

In the above calculation, $\nu$ is an arbitrary probability
measure, wherefore the equations apply to $\omega^{}_t$, for an arbitrary $t\ge 0$. 
Whenever $G\cap I \neq \varnothing$, the time-dependent
$k$-point functions are 
polynomially dependent of $\ell$-point functions with $\ell < k$. 
Consequently, they do not contain new information.
So far, we have:
\begin{prop} \label{rule-2}
   The $k$-point function $F_G^{(t)}( j^{}_1,\dots ,j^{}_k )  = 
   T^{}_{G} (\omega^{}_t) \big(\langle j^{}_1,\dots ,j^{}_k \rangle\big)$
   can only be non-\-vanishing and $($polynomially\/$)$ independent from 
   $\ell$-point functions  with $\ell < k$ if\/ 
   $G = \overline{I} = \{ \beta < j^{}_1 \} \cup \{ \beta > j^{}_k\}$.  \qed
\end{prop}
We choose this collection of $k$-point functions as our linkage equilibria.

Let us finally observe that the summation of a $k$-point function over
all possible values $x^{}_i$ of one of the specified $X_i$ 
(i.e.\ marginalization) reduces it to a $(k\!-\! 1)$-point function, 
so we have one extra (linear) relation. This means that only 
$M_i \!-\! 1$ possible values can be prescribed independently
at site $i$. On the other hand, given 
$\langle j^{}_1,\dots ,j^{}_k \rangle$, there is only one way to choose
$G$ due to Proposition \ref{rule-2}, and then there are 
$(M^{}_{j^{}_1} - 1)\cdot ... \cdot (M^{}_{j^{}_k} - 1)$ different
and independent choices to specify the actual values at the sites. 
Summing up all these possibilities results in
$$ \sum_{D\subset N}\, \prod_{i\in D} \, (M^{}_i - 1) \; = \;
   \prod_{i=0}^n \big(1 + (M^{}_i - 1)\big) \; = \;
   \prod_{i=0}^n M^{}_i \; = \; |X|\, .
$$
This means that we have singled out the right number of functions.
In view of Proposition \ref{rule-2}, for $t$ arbitrary but fixed,
they completely determine the
value of the signed measures $T^{}_{G} (\omega^{}_t)$ on all cylinder
sets. These, in turn, are closed under finite intersections and generate
the full $\sigma$-algebra of the (finite) space $X$, so all measures
$T^{}_{G} (\omega^{}_t)$, and hence also $\omega^{}_t$, are uniquely
specified, and we have achieved our goal.
An explicit example has been worked out in Section 4 of \cite{EB},
where the $k$-point functions 
$F^{(t)}_{G} ( j^{}_1,\dots ,j^{}_k )$
appear as the components of the vector $z$ of linkage disequilibria,
up to a change of basis in the local site spaces.

If $X$ is not a finite set, one has to use a generating family of Borel
cylinder sets instead of just singleton sets, and invoke Fact \ref{baire}.
Although there is no simple counting argument, the general structure is 
still similar.

At this point, one could still argue that $k$-point functions w.r.t.\
the selection of {\em sites}, as our $F^{(t)}_{G} ( j^{}_1,\dots ,j^{}_k )$,
should be replaced by proper $k$-point {\em correlation\/} functions
because these separate off all contributions of functions of lower order,
i.e.\ of $\ell$-point functions with $\ell < k$. This is just another 
application of the M\"obius inversion principle, but one where {\em all\/}
partitions (rather than only ordered ones) are needed. We provide the 
corresponding formulas in the Appendix. If one performs the necessary
calculations, one quickly realizes that our previous inclusion-exclusion
process w.r.t.\ ordered partitions of the {\em links\/} has far reaching
consequences: most of the potential correction terms simply vanish, as
a result of Lemma~\ref{zero-rule}. In particular, we obtain
\begin{theorem} \label{linkage}
   Let $S=\{j^{}_1,\ldots,j^{}_k\}$ be a set of site 
   indices, in increasing order and without gaps, and let 
   $G=\{\alpha < j^{}_1\} \cup \{\alpha > j^{}_k\}$.
   Then, the $k$-point function 
   $F^{(t)}_{G} ( j^{}_1,\dots ,j^{}_k ) = T^{}_G(\omega^{}_t)
   \big(\langle j^{}_1,\ldots,j^{}_k \rangle \big)$
   coincides with the corresponding $k$-point correlation
   function as given in Eq.~$(\ref{corr2})$ of the Appendix. 

   These functions, for all possible choices of the set $S$,
   form a polynomially independent set of linkage disequilibria.
\end{theorem}
{\sc Proof}: We apply Lemma~\ref{zero-rule} with $\nu = \omega^{}_t$.
Due to the assumption on $S$ versus $G$, the right-hand side of 
Eq.~(\ref{corr2}) boils down to the one term we already have,
because all other terms vanish. 
Propositions~\ref{rule-2} and \ref{measurement} ensure
the polynomial independence of these objects, which are
our linkage disequilibria.  \qed

\smallskip
This result does not extend to all $k$-point functions. If, for a given
$k$-point function, a non-vanishing correction term occurs in the 
corresponding correlation function, this will, in general, not decay with 
the same exponential rate as the original $k$-point function. So, grouping
according to decay rates and according to correlation structures 
simultaneously is not possible in general. It is a rather remarkable
fact that the set of linkage disequlibria is a set of exceptions,
and one (as we demonstrated above for the case of discrete state spaces)
that completely determines the probability measure.

\bigskip
\bigskip

\section{Mutation and recombination} \label{both}

In this section, we will just  combine the results of the
previous two sections. This is possible because, as we will see, mutation
and recombination are independent in our approach, i.e.\ the corresponding
operators in the IPL equation commute. This is to be expected given the
fact that mutation acts on the {\em sites\/} while recombination works
via the {\em links\/}. However, to be able to formulate
this in a more general situation than $X$ finite or discrete, we now
restrict ourselves to the Banach space 
$\cM^{\otimes} = \bigotimes_{i\in N} \cM(X_i)$
which, as explained earlier, is meant as the completion of the algebraic
tensor product. In general, it is a (true) Banach subspace of $\cM(X)$.
Our IPL equation now reads
\begin{equation} \label{b-master}
   \dot{\omega} \; = \; \Bigl(\,\sum_{i\in N} \mu^{}_i \, Q^{}_{i}
   \, + \, \sum_{\alpha\in L} \varrho^{}_{\alpha} 
   \bigl(R^{}_{\alpha}- \bs{1} \bigr)\Bigr) \, (\omega)
\end{equation}
where we have taken the liberty to introduce mutation rates $\mu^{}_i$,
all of which are assumed to be strictly positive. The idea behind this is to
use some standardized version for the mutation operators $Q^{}_i$ of
(\ref{qi}) so
that the $\mu^{}_i$ serve as relative coefficients, in line with the
usual practice in the biological literature. The linear operators $Q^{}_i$
are supposed to be bounded, hence continuous, and thus possess a unique
extension to $\cM^{\otimes}$, compare \cite[Thms.~II.1.2 and II.1.5]{Werner}. 
To show consistency, we observe
\begin{lemma}\label{invar}
  The Banach space $\cM^{\otimes}$ is invariant under $R^{}_{\alpha}$,
  for all $\alpha \in L$, and hence positive invariant 
  under the flow of\/ $(\ref{b-master})$.
\end{lemma}

{\sc Proof:} It is clear that $R^{}_{\alpha}$ maps a finite linear
combination of product measures onto another linear combination of
this kind, compare the proof of Prop.\ \ref{conv-linear}. Since such
linear combinations are dense in $\cM^{\otimes}$ and $R^{}_{\alpha}$ is
Lipschitz on $\cM(X)$, it maps the closed subspace $\cM^{\otimes}$ of
$\cM(X)$ into itself.
The statement on positive invariance is a direct consequence
of \cite[Thm.~16.5 and Remark 16.6]{Amann}.  \qed

Let $\cP^{\otimes}$ be the subspace of probability measures in $\cM^{\otimes}$.
Referring back to Propositions
\ref{markov-sol} and \ref{reco-dgl} and to Theorem \ref{general-ODE}, the 
following result is immediate.
\begin{prop} \label{unique}
   The abstract Cauchy problem of the IPL equation $(\ref{b-master})$,
   with initial condition\/ $\omega^{}_0\in\cM^{\otimes}$, has a unique solution.
   The cone $\cM^{\otimes}_{+}$ is positive invariant, and the norm of a positive
   measure is preserved in forward time. In particular, the convex set\/
   $\cP^{\otimes}$ is positive invariant.  \qed
\end{prop}

To continue, let us call a positive linear operator $W$ on $\cM^{\otimes}$
{\em strictly positive\/} if $\omega\in\cM^{\otimes}_{+}$ with $\omega>0$ implies
$W\omega > 0$. The key observation is now
\begin{lemma} \label{all-commute}
   Let $W$ be a strictly positive bounded linear operator on $\cM^{\otimes}$
   which has a complete tensor product structure, i.e.\ 
   $W=w^{}_0\otimes \dots\otimes w^{}_n$.
   On $\cM^{\otimes}_{+}$, the elementary recombinator $R^{}_{\alpha}$ then commutes
   with $W\!$, i.e.\ $W R^{}_{\alpha} = R^{}_{\alpha} W$. In particular,
   this is true if\/ $W = \exp(tQ^{}_i)$ is an element of a Markov semigroup,
   as in Section $\ref{mut}$, for any $t\ge 0$, $i\in N$ and $\alpha\in L$.
\end{lemma}
{\sc Proof}: Let us first consider the case that $W$ preserves the norm
of a positive measure $\nu$, i.e.\ $\|W\nu\| = \|\nu\|$. This is also true
of $R^{}_{\alpha}$, $\alpha\in L$. Since $W$ is linear and $R^{}_{\alpha}$
positive homogeneous of degree 1, it is sufficient to prove the claim
on $\cP^{\otimes}$. So, let $\nu\in\cP^{\otimes}$. 
$W$ has a complete tensor product structure, 
so $W = W^{}_{\s\alpha}\otimes W^{}_{\g\alpha}$ in particular.
Observe first that
$W^{}_{\s\alpha}\circ\pi^{}_{\s\alpha}=\pi^{}_{\s\alpha}\circ W$
and $W^{}_{\g\alpha}\circ\pi^{}_{\g\alpha}=\pi^{}_{\g\alpha}\circ W$.
These relations certainly hold when applied to a product measure
$\nu = \nu^{}_{\s\alpha}\otimes\nu^{}_{\g\alpha}$, but, due to linearity
of all mappings involved here, also on arbitrary (finite) linear
combinations of measures of this kind. The latter are dense in 
$\cM^{\otimes}$, so that continuity of the mappings establishes
the relations, compare \cite[Thm.~II.1.5]{Werner}.

As a consequence, we obtain
$$ W \big( (\pi^{}_{\s\alpha}.\nu)\otimes (\pi^{}_{\g\alpha}.\nu)\big) 
   \; = \; \big( \pi^{}_{\s\alpha}.(W \nu)\big) \otimes
           \big( \pi^{}_{\g\alpha}.(W \nu)\big)\,, $$
which proves the assertion for the case that $W$ preserves the norm of $\nu$.

Let us now consider the general case. The proof so far only required
that $W$ preserved the norm of the single $\nu$ under consideration.
We employ again positive homogeneity of $R^{}_{\alpha}$.
If $\nu>0$, we have $W\nu>0$ by assumption, so that 
$a:=\frac{\|\nu\|}{\|W\nu\|} > 0$ is well defined. So we obtain
$\|\nu\| = \|a W \nu \|$ and
\begin{eqnarray*}
   W R^{}_{\alpha}(\nu) & = &
   \mbox{$\frac{1}{a}$}\, (a W ) R^{}_{\alpha}(\nu)  \\
   & = &  \mbox{$\frac{1}{a}$}\, R^{}_{\alpha} ( a W \nu )
          \qquad \mbox{(by above argument)}  \\
   & = &  R^{}_{\alpha} (W \nu)
          \qquad \quad \; \mbox{(by Eq.~(\ref{pos-hom}))}
\end{eqnarray*}
which proves the first assertion.

The second claim is obvious because elements of a Markov semigroup are
strictly positive and because the generators $Q_i$, compare Eq.~(\ref{qi}),
have the required product structure.  \qed

\smallskip
We can now put together our previous efforts. The obvious form of the
solution of (\ref{b-master}) is now
\begin{equation} \label{ansatz2a}
   \omega^{}_t \; = \; \exp\bigl( t Q \bigr) \,
   \sum_{G\subset L} \, a^{}_{G}(t)\, R^{}_{G} \, (\omega^{}_0)
\end{equation}
with $Q=\sum_{i=0}^{n}\, \mu^{}_i Q^{}_i$ and
the coefficient functions $a^{}_{G}(t)$ of Theorem \ref{sol1}.
The verification that this indeed solves the IPL equation is a
simple application of the product rule. Let 
$\nu^{}_t := \sum_{G\subset L} \, a^{}_{G}(t)\, R^{}_{G} (\omega^{}_0)$,
so that $\omega^{}_t = \exp(tQ)\, \nu^{}_t$. Then we have
\begin{eqnarray*}
   \dot{\omega}^{}_t  & = & 
   Q \,\omega^{}_t + \exp\bigl( t Q \bigr) \dot{\nu}^{}_t \\[2mm]
   & = & Q \,\omega^{}_t + \exp\bigl( t Q \bigr)
     \sum_{\alpha\in L} \varrho^{}_{\alpha} \bigl(R^{}_{\alpha}-\bs{1}\bigr)
     (\nu^{}_t) \qquad \quad  \mbox{(by Theorem \ref{sol1})} \\
   & = & Q \, \omega^{}_t + \Bigl(\,\sum_{\alpha\in L} 
   \varrho^{}_{\alpha} \bigl(R^{}_{\alpha}-\bs{1}\bigr)\Bigr)
   \bigl(\exp(t Q ) \nu^{}_t \bigl)
   \qquad \, \mbox{(by Lemma \ref{all-commute})} \\
   & = & \Bigl(\,\sum_{i\in N} \mu^{}_i \, Q^{}_{i}
   \, + \, \sum_{\alpha\in L} \varrho^{}_{\alpha} 
   \bigl(R^{}_{\alpha}-\bs{1}\bigr)\Bigr) \, (\omega^{}_t)\, .
\end{eqnarray*}
So, together with Proposition \ref{unique}, we have established:
\begin{theorem} \label{full-sol}
   The unique solution of the IPL equation $(\ref{b-master})$, with
   initial condition\/ $\omega^{}_0\in\cM^{\otimes}_{+}$, is given by\/ $\omega^{}_t$
   of\/ $(\ref{ansatz2a})$, with the coefficient functions\/ $a^{}_{G}(t)$ of 
   Theorem $\ref{sol1}$. \qed
\end{theorem}

Let us take a closer look at the asymptotic behaviour. Since $a^{}_G(t)$
decreases, as $t\to\infty$, exponentially to $0$ unless $G=L$, we obtain
\begin{eqnarray*}
  \omega^{}_t & \sim & \exp(tQ) \big( a^{}_L(t) R^{}_L (\omega^{}_0)\big)
    \; \sim \; \exp(tQ) \, \bigotimes_{i=0}^{n}\, (\pi^{}_i . \omega^{}_0) \\
    & = & \bigotimes_{i=0}^{n}\, \big( \exp(t q^{}_i)
    (\pi^{}_i . \omega^{}_0) \big)  
\end{eqnarray*}
where all neglected terms are of lower order in that they vanish exponentially 
(recall that $\exp(tQ)$ and $\exp(t q^{}_i)$ are Markov). This shows that the
stationary measure, for any initial measure $\omega^{}_0$, is again a complete 
product measure\footnote{Convergence to product measures is also known from 
various interacting particle systems, compare \cite{BAZ}.}.
Whether or not there is a unique global equilibrium measure then depends on the
properties of the local mutation operators $q^{}_i$. In the case that
$X$ is finite, uniqueness follows if all these generators are irreducible.

What remains to be done is to extend the M\"obius trick and to evaluate
the linkage disequilibria also for this case. Due to Lemma \ref{all-commute}, 
we can equivalently write $\omega^{}_t$ of (\ref{ansatz2a}) as
\begin{equation} \label{ansatz2b}
   \omega^{}_t \; = \; 
   \sum_{G\subset L} \, a^{}_{G}(t)\, R^{}_{G} \, 
   \bigl(\exp(t Q ) \,\omega^{}_0 \bigr) \, .
\end{equation}
At any fixed instant of time, $\exp(tQ)\, \omega^{}_0$ is a positive measure, 
and we can employ Eq.~(\ref{new-form}) to obtain
\begin{equation} \label{dis-combined}
   \omega^{}_t \; = \; 
   \sum_{K\subset L} \, b^{}_{K}(t)\, T^{}_{K} \, 
   \bigl(\exp ( t Q ) \,\omega^{}_0 \bigr) 
\end{equation}
with the functions $b^{}_K(t)$ introduced in (\ref{new-coeff}).

If we now assume again that $X$ is finite, we can use the $k$-point
cylinder functions as before to select a finite set of linkage
disequilibria that completely determine the solution $\omega^{}_t$.
They are the functions
\begin{equation}
   F_G^t(j^{}_1,\ldots,j^{}_k) \; = \; T^{}_G (\omega^{}_t)
    \big(\langle j^{}_1, \dots , j^{}_k \rangle \big)
\end{equation}
for $G\subset L$ and selected cylinder sets 
$\langle j^{}_1, \dots , j^{}_k \rangle$ exactly as before.

Since mutation and recombination are independent of each other and
the time evolutions commute, we can separate the time decay due to
the two processes. The effect is as follows. Recombination is
sensitive to sites selected in the cylinder sets, but not to the
actual values prescribed there. Mutation, in turn, has a tensor
product structure with respect to the sites (which expresses
the independence of individual  events).

If $\exp(tQ)$ is Markov (so that Lemma $\ref{all-commute}$ applies),
it is easy to derive (in analogy with the proof of Proposition 
$\ref{measurement}$) that
\begin{equation} \label{TGQ}
   \frac{d}{dt} T^{}_G(\omega^{}_t)  \; = \; 
   \Big( Q - \sum_{\alpha \in \overline{G}} \varrho^{}_{\alpha} \Big)\,
   T^{}_G(\omega^{}_t)\,.
\end{equation}
This shows how the recombination rates and the eigenvalues of $Q$
together determine the fine structure of exponential decay.
Note that diagonalizing $Q$ (if at all possible) now corresponds 
to taking appropriate linear combinations of 
$T^{}_G (\omega^{}_t) \big(\langle j^{}_1, \dots , j^{}_k \rangle \big)$
for fixed $G$ and $j^{}_1, \dots , j^{}_k$, but different values
prescribed at the sites. For finite $X$, this has been worked out
in \cite{EB}, along with explicit examples.

\bigskip
\bigskip
\section{Selection} \label{sel}

Let us first look at selection in a slightly more general way, i.e.\ via
an IPL equation on $\cM(X)$ without explicit reference to its tensor product
structure. Let $P\! : \cM(X) \to \cM(X)$ 
be a bounded linear operator which generates a positive semigroup. 
According to \cite[Thm.\ 1.11]{Arendt}, the latter is true if and only 
if $P$ satisfies our assumption (A2), the positive minimum principle.
Consider now the ODE
\begin{equation} \label{s-master}
   \dot{\omega} \; = \; \Phi^{}_{\rm sel} (\omega) \; := \;
    P\omega - \frac{ P \omega (X)}{\|\omega\|}\, \omega 
\end{equation}
where $\Phi^{}_{\rm sel}(0)=0$ is the proper extension of 
$\Phi^{}_{\rm sel}(\omega)$ to $\omega=0$.
This is motivated by the standard selection model 
(cf.\ \cite{Hofbauer}), where, in properly 
coordinatized form as indicated in Section~\ref{prel},  
$P$ is a diagonal matrix which keeps track of the 
`fitness' of the various states, and $\frac{P \omega (X)}{\|\omega\|}$ 
is the `mean fitness' of the population. 
This model also arises in the infinite population limit of the
well-known Moran model, see \cite{Kingman} or \cite[Ch.\ 3]{Ewens}. 
Here, in a population of $m$ individuals with finite state space $X$
 as described in Section~\ref{recombi}, every
individual of type $x$ reproduces at rate $r_x$, and the offspring
replaces a randomly chosen individual in the population (possibly its
own parent). Therefore, a transition from population state
$z$ to $z + u_x - u_y$
occurs at rate $r_x z_x z_y / m$. Along the lines of Section 
\ref{recombi}, the limit $m \to \infty$ yields a special case of
the differential equation (\ref{s-master}), where $P$ is the diagonal 
matrix with elements $r_x$.

The more general form used
here does not only cover more general $X$, but also 
interaction between mutation and reproduction
(as opposed to the independent processes considered so far),
e.g.\ the production of mutated offspring on the occasion of reproduction.
In any case, the subtraction of the second term on the 
right hand side of (\ref{s-master})
comes from the preservation of total mass, or, in more
technical terms, is designed so that $\Phi^{}_{\rm sel}$ satisfies 
assumptions (A2) and (A3) from Section \ref{mut}.

So far, our selection equation seems to imply that selection acts
on haploids (i.e.\ individuals with only one copy of the genetic information
per cell). If, however, individuals have two copies that are
equivalent and do not interact (the diploid case without
dominance), Eq.\ (\ref{s-master}) is replaced by
\begin{equation} \label{diploid}
       \dot \omega \; = \; \frac{ 
       M(\omega \otimes P \omega + P \omega \otimes \omega)}
        {\|\omega\|} \;  - \;  \frac{  
       \big(M(\omega \otimes P \omega + P \omega \otimes \omega)\big) (X)}
       {\|\omega\|^2}  \; \omega  \, ,   
\end{equation}
where $M(\mu \otimes \nu) := \nu(X) \cdot \mu$ denotes marginalization
with respect to the second factor. In this formulation, the mean fitness is
$$   \frac{1}{\|\omega\|^2}\, 
    \big(M(\omega \otimes P \omega + P \omega \otimes \omega)\big) (X) 
    \; = \; 2\, \frac{\omega(X) P \omega (X)}{\|\omega\|^2} \, . $$
{}For positive $\omega$, the right-hand side of (\ref{diploid}) becomes
$$ \frac{P \omega(X)}{\|\omega\|} \, \omega \, + \, P \omega
   \, - \, 2\, \frac{P \omega(X)}{\|\omega\|} \, \omega 
   \; = \; P \omega - \frac{P \omega(X)}{\|\omega\|} \, \omega \, ,
$$
that is, the diploid equation reduces to the haploid one in this
case, in the sense that the flow is the same on $\cM_+(X)$.

Let us now take a closer look at the differential equation (\ref{s-master}).
\begin{fact} \label{s-lip}
    The mapping\/ $\Phi^{}_{\rm sel}\!: \cM(X) \to\cM(X)$ is\/ 
    $($globally\/$)$ Lipschitz.
\end{fact}
{\sc Proof}: Consider $\omega,\omega' \in \cM(X)$. If one of them is the zero 
measure, $\omega'$ say, we get
\begin{eqnarray*}
  \big\| \Phi^{}_{\rm sel}(\omega) - \Phi^{}_{\rm sel}(0) \big\| & = &
  \big\| \Phi^{}_{\rm sel}(\omega) \big\| \; \le \;
  \| P\omega \| + \frac{|P\omega (X)|}{\|\omega\|}\, \|\omega\|  \\
  & \le & 2\, \| P\omega \| \; \le \; 2 \, \|P\| \, \|\omega\|
\end{eqnarray*}
where $\|P\| := \sup_{\|\omega\|\le 1} \|P\omega\| < \infty$ because
$P$ is a bounded operator by assumption, and clearly
$|P\omega (X)|\le |P\omega|(X) = \|P\omega\|$.

Let now $\omega,\omega'$ both be non-zero. Then
$$ \big\| \Phi^{}_{\rm sel}(\omega) - \Phi^{}_{\rm sel}(\omega') \big\|
   \; \le \; \| P \| \, \| \omega - \omega' \| +
   \left\| \frac{P\omega' (X)}{\|\omega'\|}\, \omega' -
           \frac{P\omega (X)}{\|\omega\|}\, \omega \right\|. $$
Observe that $\frac{P\omega (X)}{\|\omega\|} = 
P\big(\frac{\omega}{\|\omega\|}\big) (X)$. The second term
on the right hand side of the above equation is then clearly majorized by
$ \|P\| \, \|\omega' - \omega \| + c\, \|\omega\|$ where
\begin{eqnarray*}
   c & = & \left|\,
   P\Big(\frac{\omega'}{\|\omega'\|}\Big) (X) -
   P\Big(\frac{\omega}{\|\omega\|}\Big) (X)
   \, \right| \; = \; \left| P \Big(
   \frac{\omega'}{\|\omega'\|} - \frac{\omega}{\|\omega\|}
   \Big) (X) \, \right|     \\[2mm]    & \le & \left\| P \Big(
   \frac{\omega'}{\|\omega'\|} - \frac{\omega}{\|\omega\|}
   \Big) \right\| 
   \;\, \le \;\, \frac{\|P\|}{\|\omega\|\,\|\omega'\|}\;
   \Big\|\, \|\omega\|\, \omega' - \|\omega'\|\, \omega \, \Big\|\, .
\end{eqnarray*}
Next, observe that
$$ \Big\|\, \|\omega\|\, \omega' - \|\omega'\|\, \omega \, \Big\|
   \; \le \; \|\omega'\|\, \big|\|\omega\| - \|\omega'\|\big|
   + \|\omega'\|\,\|\omega'-\omega\| \; \le \;
   2\, \|\omega'\|\,\|\omega - \omega'\| $$
so that we finally get
$$ \big\| \Phi^{}_{\rm sel}(\omega) - \Phi^{}_{\rm sel}(\omega') \big\|
   \; \le \; 4 \, \|P\| \, \| \omega - \omega'\| \, .$$
Together with the previous calculation, we see that $\Phi^{}_{\rm sel}$
is globally Lipschitz, with Lipschitz constant $\le 4$. \qed

\smallskip
So, we know that the IPL equation (\ref{s-master}) defines a unique flow.
As before, we have to check what happens with $\cM_+(X)$ under the semiflow
in forward time. Since $\omega^{}_0=0$ trivially implies $\omega^{}_t\equiv 0$
for all $t\ge 0$, we exclude this case from now on. Note that $\omega^{}_0\neq 0$
results in $\|\omega^{}_t\| > 0$ for all $t\ge 0$, due to uniqueness.
Let $\omega\in\cM_+(X)$ be a positive measure and $E$ a 
Borel set such that $\omega(E)=0$. This implies 
$\Phi^{}_{\rm sel} (\omega)\,(E) = (P\omega)\,(E) \ge 0$ 
because $P$ itself satisfies the positive minimum
principle (A2) by assumption. Also, for any $\omega\in\cM_+(X)$, we have
$$ \Phi^{}_{\rm sel} (\omega)\,(X) \; = \; P\omega (X) -
   \frac{P\omega (X)}{\|\omega\|}\, \omega(X) \; = \; 0 $$
because $\omega(X)=\|\omega\|$ for positive measures.
Together with Fact \ref{s-lip}, we see that assumptions (A1) -- (A3) are 
satisfied, and we can invoke Theorem \ref{general-ODE}. 
\begin{prop} \label{s-sol}
   Assume that the linear operator $P$ is bounded and satisfies {\rm (A2)}.
   Then the abstract Cauchy problem of the IPL equation $(\ref{s-master})$
   with initial condition $\omega^{}_0$ has a unique solution. The cone
   of positive measures is positive invariant under the flow, and the norm
   of positive measures is preserved. In particular, $\cP(X)$ is
   positive invariant.    \qed
\end{prop}

\noindent {\sc Remark}:
We would like to mention that the assumption of bounded $P$
is somewhat restricted. For non-compact $X$, many interesting
selection models lead to unbounded $P$. For mutation and
selection alone, the more general situation has been investigated
in \cite{Eshel} and, more recently, in \cite{BuBo}, in the 
framework of analytic semigroups, compare \cite[Ch.\ II.4.a]{Engel}.
Our emphasis here is on the basic structure that emerges
from the interaction with recombination; this will also
carry over to more general cases.

Before we proceed, let us make the following observation.
\begin{fact} \label{shifts}
   If the linear operator\/ $P$ is bounded and satisfies the positive
   minimum principle, the same is true of\/ $P' = P + c\, \bs{1}$ for
   arbitrary\/ $c\in\RR$. Furthermore, the flow of the IPL equation
   $(\ref{s-master})$ on $\cM_+(X)$ remains unchanged if $P$ is
   replaced by $P'$.
\end{fact}
{\sc Proof}: If $\nu$ is a positive measure and $E$ a Borel set with
$\nu(E)=0$, then $P'\nu(E)= P\nu(E) + c\,\nu(E)=P\nu(E)\ge 0$ because
$P$ satisfies (A2) by assumption. Since $P'$ is still bounded, the IPL
equation (\ref{s-master}) with $P'$ in place of $P$ conforms to
Proposition \ref{s-sol}. If $\omega\in\cM_+(X)$, we obtain
$$ P'\omega - \frac{P'\omega (X)}{\|\omega\|}\, \omega \; = \;
   P\omega + c\, \omega - \frac{P\omega (X)}{\|\omega\|}\, \omega 
   - \frac{c\, \omega(X)}{\|\omega\|}\, \omega \; = \;
   P\omega - \frac{P\omega (X)}{\|\omega\|}\, \omega $$
from which the claim follows.  \qed

\smallskip
Once again, although the ODE (\ref{s-master}) is nonlinear, it can be solved
in closed terms. This time, we employ Thompson's trick \cite{Colin} through
the substitution
\begin{equation} \label{thompson}
   \eta^{}_t \; = \; \vartheta(t)\, \omega^{}_t
\end{equation}
where $\omega^{}_t$ is a solution of (\ref{s-master}). One then obtains
$$ \dot{\eta}^{}_t \; = \; \Big( \dot{\vartheta}(t) - \vartheta(t) 
   \frac{P\omega^{}_t (X)}{\|\omega^{}_t\|} \Big) \omega^{}_t + P \eta^{}_t $$
and a significant simplification is reached if the term in brackets vanishes
because the remaining ODE is then linear. This is achieved by the choice
\begin{equation} \label{thompson2}
   \vartheta(t) \; = \; \exp \Bigg ( \int_0^t 
   \frac{P\omega^{}_{\tau}(X)}{\|\omega^{}_{\tau}\|} \, {\rm d}\tau \Bigg )
   \; = \;  \exp \, \Bigg ( \frac{1}{\|\omega^{}_0\|} \int_0^t 
    P\omega^{}_{\tau} (X) \, {\rm d}\tau \Bigg )
\end{equation}
where the second step follows from Proposition \ref{s-sol}. Clearly,
$\vartheta(t)$ is well defined (whenever $\omega^{}_0\neq 0$, which is all
we need), and we have reduced the Cauchy problem of (\ref{s-master}) to 
that of the simple linear evolution equation
\begin{equation} \label{s-master-lin}
   \dot{ \eta }  \; = \; P \eta \, .
\end{equation}
This ODE defines a uniformly continuous positive semigroup (since $P$ was 
assumed to be bounded and to satisfy (A2), the positive minimum principle). 
The solution of (\ref{s-master-lin}) will no longer have
fixed norm, but one can always get back to $\omega^{}_t$ via
$$ \omega^{}_t \; = \; \frac{\|\eta^{}_0\|}{\|\eta^{}_t\|}\, \eta^{}_t \, .$$
Note that $\eta^{}_0 = \omega^{}_0$ and $\|\omega^{}_t\|\equiv\|\eta^{}_0\|$.

Let us next consider the function
\begin{equation} \label{lapu1}
   L(t) \; = \; \frac{P\omega^{}_t (X)}{\|\omega^{}_t\|}
\end{equation}
which is defined on any orbit of the flow of (\ref{s-master}).
$L(t)$ is of particular interest on orbits of positive measures,
where it admits the interpretation as {\em mean}\/ (or {\em averaged\/})
{\em fitness\/}. Here, we know $\|\omega^{}_t\|\equiv \|\omega^{}_0\|$ 
by Proposition \ref{s-sol}, so that we obtain
\begin{eqnarray*}
   \dot{L}(t) & = & \frac{\frac{\rm d}{{\rm d}t} 
   \big( P\omega^{}_t (X) \big)}{\|\omega^{}_t\|}
   \; = \; \frac{P\dot{\omega}^{}_t (X)}{\|\omega^{}_t\|} 
   \; = \; \frac{1}{\|\omega^{}_t\|}\, P \Big( P\omega^{}_t
   - \frac{P\omega^{}_t (X)}{\|\omega^{}_t\|}\, \omega^{}_t  
   \Big) (X)    \\[2mm]
   & = & P_{}^2 \Big( \frac{\omega^{}_t}{\|\omega^{}_t\|} \Big) (X)
    -  \Big(  P \Big( \frac{\omega^{}_t}{\|\omega^{}_t\|} \Big) (X)
    \Big)^2
\end{eqnarray*}
which has the form of a variance. So we can state
\begin{prop} \label{lapu-prop}
   If, under the assumptions of Proposition $\ref{s-sol}$, 
   $P$ satisfies the condition $(P\omega(X))^2\le P^2\omega(X)$
   on $\cP(X)$, the function $L(t)$ of $(\ref{lapu1})$ is a Lyapunov
   function for the flow of\/ $(\ref{s-master})$ on the positive cone
   $\cM_+(X)$.
\end{prop}
{\sc Proof}: From the above calculation, it is clear that
$\dot{L}(t) \ge 0$ on all orbits in $\cM_+(X)$ if $P$ satisfies
the inequality $(P\omega(X))^2\le P^2\omega(X)$ on $\cP(X)$,
so $L$ cannot decrease along such an orbit.  \qed

\smallskip
\noindent
{\sc Remark}: Our definition of a Lyapunov function on $\cM_+(X)$ is
{\em global\/} and (up to a sign) that
of \cite[Ch.\ 18]{Amann}. Note that the stricter version of \cite{GH},
where $\dot{L}(t) = 0$ would correspond to a unique equilibrium on
$\cM_+(X)$, is not so useful here because the asymptotic state 
(as $t\to\infty$) of the selection equation depends on the initial condition, 
i.e.\ there is no unique equilibrium in general. However, one might profit
from the use of {\em local\/} Lyapunov functions, compare
\cite[Thm.\ 1.0.2 (iii)]{GH}, but we do not expand on this here.

The condition on $P$ can be reformulated by noting that
$$ P^2\omega(X) - (P\omega(X))^2 \; = \; (P-c\,\bs{1})^2 \omega (X)$$
with $c = P\omega(X)$. A sufficient condition for Proposition \ref{lapu-prop}
to hold is then that $(P-c\,\bs{1})^2$ is a positive operator for all
(or sufficiently many) $c\in\RR$. A particularly well studied case of this
is when $X$ is finite and $P$ is a diagonal matrix in the canonical basis
consisting of the extremal measures of $\cP(X)$. In this case, Proposition
\ref{lapu-prop} is known as Fisher's fundamental theorem, see
\cite{Hofbauer} for details. In the more general case, Lyapunov
functions may be considered even more important since they
determine the `direction' of the evolution process in a situation
where little information is available otherwise,
since the solution given by Eq.\ (\ref{selec-sol}) is not very explicit
then.

The results of this section can also be formulated for the (sub-)space
$\cM^{\otimes}$ of $\cM(X)$, if it is invariant under the action of $P$.
In view of the product structure of $X$, let us now assume that we have
$P=\sum_{i=0}^{n} P_i$, with bounded $P_i$ that are locally represented 
by $p^{}_i$ (as with $q^{}_i$ versus $Q_i$ before). Clearly, $P$ maps 
$\cM^{\otimes}$ into itself. We call this situation additivity across 
sites, in complete analogy to our previous discussion of mutation.
We can then rewrite our solution as
\begin{equation} \label{selec-sol}
   \eta^{}_t \; = \; \exp(t P)\, \eta^{}_0 \; = \;
   \Big( \bigotimes_{i=0}^{n} \exp(t p^{}_i) \Big) \eta^{}_0 \, .
\end{equation}
 
With some further restrictions on the linear operator $P$, an analogue of
Proposition \ref{lapu-prop} remains true even in the presence of recombination.
This rests on the applicability of Lemma \ref{all-commute}. We thus consider
the IPL equation
\begin{equation} \label{mix-master}
   \dot{\omega} \; = \; P\omega - \frac{P\omega (X)}{\|\omega\|}\, \omega
   + \sum_{\alpha\in L} \varrho^{}_{\alpha} (R^{}_{\alpha} - \bs{1}) (\omega)
   \; = \; \Phi^{}_{\rm sel}(\omega) + \Phi^{}_{\rm rec}(\omega)
\end{equation}
whose Cauchy problem has all the nice properties we need, see Proposition
\ref{gen-Cauchy} below in the special case $Q=0$. We now assume:
\begin{enumerate}
\item $P$ has complete product structure as a generator, i.e.\ 
      $P = \sum_{i=0}^{n} P_i$ with
      $P_i = \bs{1} \otimes \dots \otimes \bs{1} \otimes p^{}_i
             \otimes \bs{1} \otimes \dots \otimes \bs{1}$.
\item Each $P_i$ is itself a bounded, strictly positive operator.
\end{enumerate}

If $\omega^{}_t$ is a solution of (\ref{mix-master}), we again define $L(t)$
as in (\ref{lapu1}) and obtain, by Lemma \ref{all-commute},
\begin{eqnarray*}
   \dot{L}(t) & = & \frac{1}{\|\omega^{}_t\|}\, P \Big(
   P \omega^{}_t - \frac{P\omega^{}_t (X)}{\|\omega^{}_t\|}\, \omega^{}_t
   + \sum_{\alpha\in L} \varrho^{}_{\alpha} (R^{}_{\alpha} - \bs{1}) 
   (\omega^{}_t) \Big) (X)  \\   & = &
   \bigg[ P_{}^2 \Big( \frac{\omega^{}_t}{\|\omega^{}_t\|} \Big) (X)
    -  \Big(  P \Big( \frac{\omega^{}_t}{\|\omega^{}_t\|} \Big) (X)
    \Big)^2 \bigg] \, + \;  \frac{\sum_{i=0}^{n}\Phi^{}_{\rm rec}
    ( P^{}_i\, \omega^{}_t ) (X)}{\|\omega^{}_t\|} \, .
\end{eqnarray*}
The last term vanishes due to our general assumptions because $P^{}_i\, \omega^{}_t >0$
and then $\Phi^{}_{\rm rec}( P^{}_i\, \omega^{}_t ) (X) = 0$ due to (A3). So, we 
are back to the condition already encountered above. To summarize:
\begin{theorem} \label{full-lapu}
   Let\/ $P=\sum_{i=0}^{n} P_i$ satisfy the assumptions of Proposition 
   $\ref{s-sol}$, and let each $P_i$ be a bounded, strictly positive operator 
   with complete product structure. If $P$ also satisfies the condition
   $(P\omega(X))^2\le P^2\omega(X)$ on $\cP^{\otimes}$, the function $L(t)$ 
   of $(\ref{lapu1})$ is a Lyapunov function for the flow of\/ 
   $(\ref{mix-master})$ on the positive cone $\cM^{\otimes}_{+}$.  \qed
\end{theorem}

\smallskip
In the absence of recombination, there are other Lyapunov functions known 
for certain combinations of selection with {\em mutation}. They rely on the 
spectral theorem applied to $P+Q$, see \cite{Jones}. In 
selection-{\em recombination}
equations where $P$ violates the product structure, the mean fitness $L(t)$
need no longer be a Lyapunov function. Moreover, the possibility of periodic
solutions \cite{Akin} demonstrates that, in more general
(diploid) models (e.g.\ with dominance), no meaningful 
Lyapunov function is to be expected.

\bigskip
\bigskip
\section{All three} \label{all}

In this last step, we combine all three processes, with the general 
assumptions as before. In view of the inherent product structure, we
only consider the dynamics on the Banach space $\cM^{\otimes}$.
The IPL equation now reads
\begin{eqnarray} \label{all-master}
   \dot{\omega} & = & 
   \Phi^{}_{\rm mut} (\omega) + \Phi^{}_{\rm rec} (\omega) 
    + \Phi^{}_{\rm sel} (\omega) \\[1mm]
   & = &   ( Q + P ) \omega - \frac{P \omega (X)}{\|\omega\|}\, \omega +
   \sum_{\alpha\in L} \varrho^{}_{\alpha} \big(R^{}_{\alpha} - \bs{1}\big) 
   (\omega)   \nonumber
\end{eqnarray}
and we immediately get the following result, again from Theorem
\ref{general-ODE}, and Lemma $\ref{invar}$.
\begin{prop} \label{gen-Cauchy}
   Let  $Q$ be a bounded Markov generator and  $P$ a bounded generator
   of a positive semigroup, both of product form. Let
   $R^{}_{\alpha}$ be the recombinators of Eq.~$(\ref{reco-op1})$. Then,
   the abstract Cauchy problem of the IPL equation $(\ref{all-master})$
   has a unique solution. The cone $\cM^{\otimes}_{+}$ is positive  
   invariant and
   the norm of positive measures is preserved under the forward flow. \qed
\end{prop}

\noindent
{\sc Remark}: Since $Q$ is a Markov generator, we know from Section \ref{mut}
that $Q\omega(X)=0$ for all $\omega\in\cM^{\otimes}_{+}$, and we could also 
start from an IPL equation where $Q$ is absorbed into $P$ --- it would give 
the same flow on $\cM^{\otimes}_{+}$. We retain the separation into mutation 
and selection because, in more general situations,
it is often adequate from both the biological and the mathematical
point of view (for example, the mutation operator is usually bounded,
but the selection operator may be unbounded);  for review, see
\cite[Ch.\ IV]{Reinhard}. We will also combine $Q$ and $P$, but
only after Thompson's linearization transformation.

Let $\omega^{}_t$, $t\ge 0$, be the solution for initial condition $\omega^{}_0$.
Define $\eta^{}_t$ as above in (\ref{thompson}), with $\vartheta(t)$ of 
(\ref{thompson2}). Then, $\omega^{}_t$ is a solution of (\ref{all-master}) if
and only if $\eta^{}_t$ solves the reduced IPL equation
\begin{equation} \label{all-master-red}
   \dot{\eta} \; = \; S \eta + \Phi^{}_{\rm rec}(\eta)
\end{equation}
where $S=Q+P$ is the bounded generator of a uniformly continuous semigroup of
positive operators. Note that the right hand side of (\ref{all-master-red})
still satisfies assumptions (A1) and (A2), but no longer (A3). So, the
corresponding Cauchy problem still has a unique solution, with $\cM^{\otimes}_{+}$
being positive invariant, but the norm of positive measures need no longer
be preserved under the flow in forward time --- and this is precisely the
point of this exercise!

{}From now on, we generally assume that both mutation and selection 
are adapted to the special product form of our state space $X$, so 
$S=\sum_{i=0}^n S_i$ (with corresponding local operator $s^{}_i$). 
Hence, $\exp(tS)$ is again a tensor product of local operators.
\begin{lemma} \label{full-commute}
   If\/ $S=\sum_{i=0}^n S_i$ is the bounded generator of a uniformly continuous
   semigroup of positive operators, then we have\/
   $\exp(tS) R^{}_{\alpha} = R^{}_{\alpha} \exp(tS)$  on $\cM^{\otimes}_{+}$,
   for all\/ $t\ge 0$ and $\alpha\in L$.
\end{lemma}
{\sc Proof}: Fix $t\ge 0$ and set $W=\exp(tS)$. This is a positive
operator by assumption. Also, since $S$ is bounded, $\nu>0$ implies 
$\exp(tS)\nu > 0$ and $W$ is strictly positive. The result then follows
from Lemma \ref{all-commute}.  \qed

\smallskip
This result means that we can use all our above methods again and construct
immediately the solution of (\ref{all-master-red}). At this point, we
particularly profit from our approach in that we can still solve the 
case with (additive) selection. In the context of Haldane linearization,
any form of selection  has, so far, appeared as a major obstacle, due to 
the fact that the flow induced by $P$ fails to preserve
the norm of positive measures \cite{Ringwood}.

\begin{theorem} \label{the-end}
   If\/ $S=\sum_{i=0}^n S_i$ satisfies the assumptions of 
   Lemma~$\ref{full-commute}$, the solution of the reduced
   IPL equation $(\ref{all-master-red})$, with initial condition\/
   $\eta^{}_0\in\cM^{\otimes}_{+}$, is given by
$$ \eta^{}_t \; = \; \exp(tS) \sum_{G\subset L} a^{}_{G}(t)\,
   R^{}_{G}(\eta^{}_0) $$
   with the coefficients\/ $a^{}_{G}(t)$ of\/ $(\ref{coeff-fun})$.
   The solution of the abstract Cauchy problem for the original
   IPL equation $(\ref{all-master})$ emerges from here via
$$ \omega^{}_t \; = \; \frac{\|\omega^{}_0\|}{\|\eta^{}_t \|}\, \eta^{}_t $$
   where\/ $\omega^{}_0 = \eta^{}_0$. If $\omega^{}_0\in\cP(X)$,
   then $\{\omega^{}_t \mid t\ge 0\}$ is a
   one-parameter family of probability measures. \qed
\end{theorem}

In line with our previous reasoning, we can determine the asymptotic behaviour,
$$ \eta^{}_t \; \sim \; \bigotimes_{i=0}^{n}\, 
   \big( \exp(t s^{}_i) (\pi^{}_i . \eta^{}_0) \big), $$
where we have used the product structure of $\exp(tS)$ and the
fact that all neglected terms, as $t\to\infty$, are exponentially small 
in comparison. The meaning for $\omega^{}_t$ is, once again, that
stationary measures are complete product measures, and the properties
of the linear operators $s^{}_i$ determine whether there is a unique global
equilibrium measure. This is connected to the general Perron-Frobenius
theory of positive operators which is rather involved in general, see
\cite[Ch.\ V.5]{Sch} and \cite{Greiner}. If, however, $X$ is finite
(so that $\cM^{\otimes}$ is finite-dimensional, and $\cM^{\otimes}=\cM(X)$) 
and all $s^{}_i$ are irreducible, there are unique $\nu^{}_i\in\cP(X^{}_i)$ 
so that $\exp(t s^{}_i)\nu^{}_i = \exp(t\lambda^{}_i)\nu^{}_i$ with 
$\lambda^{}_i\in\RR$ being the largest eigenvalue of $s^{}_i$.
In this case, as a simple calculation shows, we obtain
$$  \omega^{}_t \; \longrightarrow \; 
    \nu^{}_0\otimes\dots\otimes\nu^{}_n $$
in the $\|.\|$-topology, as $t\to\infty$, for any initial condition
$\omega^{}_0\in\cP(X)$.

Also, the following observation results immediately from Theorem \ref{the-end}.
\begin{coro}
    If an initial condition\/ $\omega^{}_0\in\cP^{\otimes}$ is a product 
    measure at link\/ $\alpha\in L$, this is also true of the
    corresponding solution\/ $\omega^{}_t$ of\/ $(\ref{all-master})$, for
    all $t\ge 0$. In particular, if\/ $\omega^{}_0$ is a complete
    product measure, this remains the case under the forward flow,
    i.e.\ for all\/ $\omega^{}_t$ with $t\ge 0$. \qed
\end{coro}

Let us return to the general discussion.
The remainder is then a copy of what we did in Section \ref{both}, 
with $Q$ replaced by $S$. In particular, we get
$$ \eta^{}_t \; = \; \sum_{K\subset L} b^{}_{K}(t)\,
   T^{}_{K} \big( \exp(tS) \eta^{}_0 \big) $$
from which one can, once again, determine the linkage disequilibria.
Note, however, that the meaning has changed now, because the norm of
$\eta^{}_t$ varies with time. In particular, one has to consider the 
quotient $\eta^{}_t/\|\eta^{}_t\|$, rather than $\eta^{}_t$ alone, 
to extract the correct behaviour for the linkage disequilibria
$F_G^t(j^{}_1,\ldots,j^{}_k) = T^{}_G (\omega^{}_t)
\big(\langle j^{}_1, \dots , j^{}_k \rangle \big)$. To be concrete, observe
first that 
$$  \frac{d}{dt} T^{}_G(\eta^{}_t) \; = \; 
    \Big( S - \sum_{\alpha \in \overline{G}} \varrho^{}_{\alpha} \Big) \,
    T^{}_G(\eta^{}_t)   $$
in perfect analogy with (\ref{TGQ}). Since $T^{}_G$ is positive
homogeneous of degree one (Eq.\ (\ref{t-hom})), and 
$\|\eta^{}_t\| = \eta^{}_t(X)$ for positive measures, one
obtains
$$   \frac{d}{dt} T^{}_G(\omega^{}_t)  \; = \;
     \Big( S - \frac{S \omega^{}_t(X)}{\|\omega^{}_t\|}\, \bs{1}
     - \sum_{\alpha \in \overline{G}}\, \varrho^{}_{\alpha} \Big)
     T^{}_G(\omega^{}_t)\,.  $$
Clearly, knowledge of the mean fitness, $S \omega^{}_t(X)/\|\omega^{}_t\|$,
is now required to determine the dynamics of the linkage disequilibria.

\bigskip
\bigskip
\section{Afterthoughts} \label{after}

In this article, we have constructed an explicit solution
of the single-crossover recombination model in continuous time, with
mutation and additive selection. It is quite astonishing
that such a solution should be possible at all -- after all,
it is an explicit representation of a nonlinear semigroup.
However, it is no coincidence that this works
in continuous time, rather than  in  discrete
time. Let us discuss this for recombination alone.
The  discrete-time analogue of our single-crossover model
is the so-called model with {\em complete interference} \cite{Christiansen2}:
\begin{equation}
   \omega^{}_{n+1}  \; = \; 
   \sum_{\alpha \in L} \varrho_{\alpha} R_{\alpha} (\omega^{}_n)
      + \Big (1 - \sum_{\alpha \in L} \varrho_{\alpha} \Big )\, \omega^{}_n\,.
\end{equation}
Similar as this may look to its continuous-time relative,
the probabilistic structure is quite different.
Single crossovers in continuous time imply {\em independence of links},
as expressed in the coefficient functions (\ref{coeff-fun})
and the resulting factorization property (Lemma \ref{coeff-factorization}).
In contrast, a second crossover is inhibited for the
duration of an entire generation in  discrete time,
due to interference of crossovers with each other (hence the name);
see also \cite{Manos}.
As a result, independence is lost, which makes the discrete
model inherently more difficult.

Of course, this also applies to the situation with selection.
Models of recombination and selection based on independent sites and 
finite site spaces have been thoroughly investigated in the population 
genetics literature, see \cite{Ewens,Ewens2,Karlin,KL,KiLy,Lyu}
for some key references and
\cite{Reinhard,Christiansen} for recent comprehensive reviews.
Independence of sites with respect to selection is reflected by a
tensor product structure of $P$,  may be interpreted as lack of
interaction between genes, and is known as absence of {\em epistasis\/}
in genetics. More precisely, since  the dynamical systems mostly considered
so far were in discrete time, a comparison with our setting
is more adequate at the level of the semigroup, rather than 
that of the generator. 

Two notions of independence have been used, compare 
 \cite{Karlin,KiLy,Reinhard}, 
which would translate into our setting as either  
$\exp(P) = \prod_i \exp(P_i) = \bigotimes_i \exp( p_i)$ 
(`multiplicative fitness') 
or as $\exp(P)$ replaced by $\sum_i \exp( P_i)$ (`additive fitness'). 
Previously, much emphasis has been on the effects of dominance
(i.e.\ the interaction between the two alleles
joined in a diploid genotype). This may lead to multiple
equilibria, which need not all be of product type,
and astonishing differences in the qualitative
behaviour of the multiplicative and 
additive scenario are observed, as reviewed in \cite{Karlin,Reinhard}.
However, these effects are absent if there is no dominance
(as in our model); in particular, all equilibria are then of product
type.
Thus, our simple continuous-time model might well serve
as an exactly solved reference case which also captures the qualitative
features of the corresponding models in discrete time, although no
explicit solution is available there.

Now, the logical next step would be to extend the analysis to the inclusion
of interactions between sites, which occur as soon as
selection is no longer additive across sites.
Alas, this is much more involved, and even the simplest cases
go far beyond what we have outlined above. The reason is that selection
now forces the introduction of further terms in the right hand side of the 
IPL equation so that the corresponding semigroups  no longer commute 
with recombination. Nevertheless, several situations
can be envisioned that admit at least a perturbative approach.
In line with the single-crossover assumption, an expansion
for small recombination rates would be appropriate,
in contrast to the well-known quasi-linkage-equilibrium
approach 
for large recombination rates (for review, see \cite{Reinhard}).
We hope to report on some progress in this direction soon.

\bigskip

\subsection*{Appendix: Moments versus correlations}

As mentioned above, it is often desirable to separate effects that stem from
mutual interactions of differently many ``particles'' or, as in the above
discussion, from specification at a different number of sites. For two sites,
correlation $C$ and moments $F$ are related by
$C(\{i,j\}) = F(\{i,j\}) - F(\{i\}) F(\{j\})$, where the arguments
are meant as symbolic labels. Since this is a rather general structure, we 
briefly describe its systematic treatment by means of M\"obius inversion,
also known as inclusion-exclusion principle.

Let $S=\{1,2,\dots , k\}$ be a finite set which will serve as the index set
of the particles or the specified sites, the latter through
$\langle j^{}_1,\dots , j^{}_{\ell}\rangle$.
Let ${\mathcal A}=\{A^{}_1,\dots ,A^{}_p\}$ be a {\em partition\/} of $S$, 
i.e.\ $S$ is the disjoint union of the members of $\mathcal A$. 
Unlike before, the partition need not be ordered. Let the partition 
${\mathcal B}=\{B^{}_1,\dots ,B^{}_q\}$ 
be a {\em refinement\/} of $\mathcal A$, so that
$$ A^{}_1=B^{}_{j^{}_{1,1}}\cup\dots\cup B^{}_{j^{}_{1,n^{}_1}} 
   \, , \; \dots \; , \,
   A^{}_p=B^{}_{j^{}_{p,1}}\cup\dots\cup B^{}_{j^{}_{p,n^{}_p}} $$
where $\{\{j^{}_{1,1},\dots ,j^{}_{1,n^{}_1}\},\dots,
\{j^{}_{p,1},\dots,j^{}_{p,n^{}_p}\}\}$ is a partition of
$\{1,\dots, q\}$, hence $n^{}_1 + \ldots + n^{}_p = q$.
We write ${\mathcal B} \preccurlyeq {\mathcal A}$ in this case, where
$\preccurlyeq$ defines a partial order which makes $S$ into a poset.
The corresponding M\"obius function, compare \cite[p.\ 86]{Berge}, is
given by 
\begin{eqnarray} \label{gen-part}
   \mu({\mathcal B},{\mathcal A}) & = & \prod_{i=1}^{p}\,
   (-1)^{n^{}_i -1} (n^{}_i - 1)! \\
   & = & (-1)^{p+n^{}_1+\ldots+n^{}_p} \, (n^{}_1 - 1)!
   \cdot\ldots\cdot (n^{}_p - 1)!  \nonumber
\end{eqnarray}
If $\mathcal C$ is any refinement of $\mathcal A$,
$\mu$ satisfies the formula
$$ \sum_{{\mathcal C}\preccurlyeq
         \text{\d{$\mathcal B$}}\preccurlyeq{\mathcal A}}
   \mu({\mathcal B},{\mathcal A}) \; = \;
   \begin{cases}  1 & \text{if ${\mathcal A}={\mathcal C}$} \\
                  0 & \text{otherwise.}    
    \end{cases}  $$

Let us now, for a partition ${\mathcal A} = \{A^{}_1,\dots ,A^{}_p\}$, 
introduce the function $F({\mathcal A})=F(A^{}_1)\cdot\ldots\cdot F(A^{}_p)$, 
and similarly for the correlations, $C$. These quantities are related by
\begin{equation} \label{corr1}
   F({\mathcal A}) \; := \; \sum_{{\mathcal B}\preccurlyeq{\mathcal A}}
   C({\mathcal B}) \; = \;  \sum_{{\mathcal B}\preccurlyeq{\mathcal A}}\,
   \prod_{B\in{\mathcal B}}\, C(B) 
\end{equation}
because this precisely reflects the idea to separate off contributions
from subsets of different cardinality.
The M\"obius inversion formula then gives the following formula for
the special case that ${\mathcal A} = \{ A \}$:
\begin{equation} \label{corr2}
   C(A) \; = \; \sum_{{\mathcal B}\preccurlyeq{\mathcal A}}
   F({\mathcal B})\, \mu({\mathcal B},{\mathcal A})
   \; = \; \sum_{{\mathcal B}\preccurlyeq{\mathcal A}}\, 
   (-1)^{|{\mathcal B}| - 1}
   \big( |{\mathcal B}| - 1 \big)! \,
   \prod_{i=1}^{|{\mathcal B}|}\, F(B_i) \;
\end{equation}
where $|{\mathcal B}|$ denotes the number of sets in the partition 
${\mathcal B} = \{ B^{}_1,\dots ,B^{}_{|{\mathcal B}|}\}$. The 
following example might illustrate this formula:
\begin{eqnarray*}
   C(\mbox{\small $\{1,2,3\}$}) & = & F(\mbox{\small $\{1,2,3\}$})  
     \; + \;  2\, F(\mbox{\small $\{1\}$}) F(\mbox{\small $\{2\}$}) 
                 F(\mbox{\small $\{3\}$}) \\
   & & - \;  F(\mbox{\small $\{1\}$}) F(\mbox{\small $\{2,3\}$}) 
              - F(\mbox{\small $\{2\}$}) F(\mbox{\small $\{1,3\}$})  
              - F(\mbox{\small $\{3\}$}) F(\mbox{\small $\{1,2\}$})  
\end{eqnarray*}
which is to be compared with
\begin{eqnarray*}
   F(\mbox{\small $\{1,2,3\}$}) & = & C(\mbox{\small $\{1,2,3\}$})  
     \; + \;  C(\mbox{\small $\{1\}$}) C(\mbox{\small $\{2\}$}) 
                 C(\mbox{\small $\{3\}$}) \\
   & & + \;  C(\mbox{\small $\{1\}$}) C(\mbox{\small $\{2,3\}$}) 
              + C(\mbox{\small $\{2\}$}) C(\mbox{\small $\{1,3\}$})  
              + C(\mbox{\small $\{3\}$}) C(\mbox{\small $\{1,2\}$})  
\end{eqnarray*}
according to (\ref{corr1}).
Let us finally remark that formula (\ref{corr2}) can be applied
factorwise if ${\mathcal A} = \{A^{}_1,\dots ,A^{}_p\}$ because then
$C({\mathcal A}) = C(A^{}_1)\cdot\ldots\cdot C(A^{}_p)$ by definition.

\bigskip

\subsection*{Acknowledgements}

It is our pleasure to thank Reinhard B\"urger, Hans-Otto Georgii,
Joachim Hermisson, Achim Klenke and Manfred Wolff for a number of 
clarifying discussions, and Ulrich Hermisson and Oliver Redner for
carefully reading the manuscript. M.B.\ would like to thank Robert V.\ Moody
and the Department of Mathematical Sciences of the University of Alberta
(Edmonton, Canada) for hospitality, where part of this work was done.


\bigskip
\bigskip

\begin{thebibliography}{99}
\small

\bibitem{Aigner}
M.~Aigner,
{\em Combinatorial Theory}, 
Springer, Berlin (1979); reprint (1997).

\vspace*{-1mm}
\bibitem{Akin-Buch}
E.~Akin,
{\em The Geometry of Population Genetics}, LNB 31, Springer, Berlin (1979).

\vspace*{-1mm}
\bibitem{Akin}
E.~Akin,
{\em Cycling in simple genetic systems},
J.\ Math.\ Biol.\ {\bf 13} (1982) 305--324.

\vspace*{-1mm}
\bibitem{Amann}
H.~Amann,
{\em Gew\"ohnliche Differentialgleichungen}, 2nd ed.,
de Gruyter, Berlin (1995); (older) English ed.:
{\em Ordinary Differential Equations},
de Gruyter, Berlin (1990).

\vspace*{-1mm}
\bibitem{Arendt}
W.~Arendt, 
{\em Characterization of positive semigroups on Banach lattices},
in: {\em One-parameter Semigroups of Positive Operators},
ed.\ R.\ Nagel, LNM 1184, Springer, Berlin (1986), pp.\ 247--291.

\vspace*{-1mm}
\bibitem{EB}
E.~Baake,
{\em Mutation and recombination with tight linkage},
J.\ Math.\ Biol.\ {\bf 42} (2001) 455--488.

\vspace*{-1mm}
\bibitem{BG}
E.~Baake and W.~Gabriel,
{\em Biological evolution through mutation, selection and drift:
An introductory review}, in:
{\em Annual Review of Computational Physics}, vol.\ VII, ed.\
D.\ Stauffer, World Scientific, Singapore (2000), pp.\ 203--264;
cond-mat/9907372.

\vspace*{-1mm}
\bibitem{BAZ}
G.~Ben Arous and O.~Zeitouni,
{\em Increasing propagation of chaos for mean field models},
Ann.\ Inst.\ H.\ Poincar\'e -- Prob.\ Stat.\ {\bf 35} (1999) 85--102.

\vspace*{-1mm}
\bibitem{Berb}
S.~K.~Berberian,
{\em Measure and Integration},
Macmillan, New York (1965).

\vspace*{-1mm}
\bibitem{Berge}
C.~Berge,
{\em Principles of Combinatorics},
Academic Press, New York (1971).

\vspace*{-1mm}
\bibitem{Reinhard}
R.~B{\"u}rger,
{\em The Mathematical Theory of Selection, Recombination, and Mutation},
Wiley, Chichester (2000).

\vspace*{-1mm}
\bibitem{BuBo}
R.~B{\"u}rger and I.~Bomze,
{\em Stationary distributions under mutation-selection balance:
structure and properties},
Adv.\ Appl.\ Prob.\ {\bf 28} (1996) 227--251.

\vspace*{-1mm}
\bibitem{Christiansen2}
F.~B.~Christiansen,
{\em The effect of population subdivision on multiple loci without
selection}, in: {\em Mathematical evolutionary theory}, ed.\ M.~W.~Feldman,
Princeton University Press, Princeton (1989), pp.~71--85.

\vspace*{-1mm}
\bibitem{Christiansen}
F.~B.~Christiansen,
{\em Population Genetics of Multiple Loci},
Wiley, Chichester (2000).

\vspace*{-1mm}
\bibitem{Clark}
A.~Clark et al.,
{\em Haplotype structure and population genetic inferences from
nucleotide-sequence variation in human lipoprotein lipase},
Am.\ J.\ Hum.\ Gen.\ {\bf 63} (1998), 595--612.

\vspace*{-1mm}
\bibitem{Dawson}
K.~J.~Dawson,
{\em The decay of linkage disequilibria under random union of gametes:
how to calculate Bennett's principal components},
Theor.\ Pop.\ Biol.\ {\bf 58} (2000) 1--20.

\vspace*{-1mm}
\bibitem{Dawson2}
K.~J.~Dawson,
{\em The evolution of a population under recombination: How to linearise the
dynamics},
Lin.\ Alg.\ Appl.\ {\bf 348} (2002) 115--137.

\vspace*{-1mm}
\bibitem{Dudley}
R.~M.~Dudley,
{\em Real Analysis and Probability},
Chapman and Hall, New York (1989).

\vspace*{-1mm}
\bibitem{Engel}
K.-J.~Engel and R.~Nagel,
{\em One-Parameter Semigroups for Linear Evolution Equations},
GTM 194, Springer, New York (2000).

\vspace*{-1mm}
\bibitem{Eshel}
I.~Eshel,
{\em Evolution processes with continuity of types},
Adv.\ Appl.\ Prob.\ {\bf 4} (1972) 475--507.

\vspace*{-1mm}
\bibitem{EK}
S.~N.~Ethier and T.~G.~Kurtz,
{\em Markov Processes: Characterization and Convergence},
Wiley, New York (1986).

\vspace*{-1mm}
\bibitem{Ewens}
W.~J.~Ewens,
\textit{A generalized fundamental theorem of natural selection},
Genetics {\bf 63} (1969) 531--537.

\vspace*{-1mm}
\bibitem{Ewens2}
W.~J.~Ewens,
{\em Mean fitness increases when fitnesses are additive},
Nature {\bf 221} (1969) 1076.

\vspace*{-1mm}
\bibitem{FW}
M.~I.~Freidlin and A.~D.~Wentzell,
{\em Random Perturbations of Dynamical Systems}, 2nd ed.,
Springer, New York (1998).

\vspace*{-1mm}
\bibitem{Greiner}
G.~Greiner,
{\em Spectral theory of positive semigroups on Banach lattices},
in: {\em One-parameter Semigroups of Positive Operators},
ed.\ R.\ Nagel, LNM 1184, Springer, Berlin (1986), pp.\ 292--332.

\vspace*{-1mm}
\bibitem{GH}
J.~Guckenheimer and Ph.~Holmes,
{\em Nonlinear Oscillations, Dynamical Systems, and Bifurcations of
Vector Fields}, corr.\ 3rd printing, Springer, New York (1990).

\vspace*{-1mm}
\bibitem{Hofbauer}
J.~Hofbauer,
{\em The selection-mutation equation},
J.\ Math.\ Biol.\ {\bf 23} (1985) 41--53.

\vspace*{-1mm}
\bibitem{Jones}
B.~L.~Jones,
{\em Some principles governing selection in self-reproducing macromolecular
systems -- an analog of Fisher's fundamental theorem},
J.\ Math.\ Biol.\ {\bf 6} (1978) 169--175.

\vspace*{-1mm}
\bibitem{vankamp}
N.~G.~van Kampen,
\textit{Stochastic Processes in Physics and Chemistry},
North-Holland, Amsterdam (1981).

\vspace*{-1mm}
\bibitem{Karlin}
S.~Karlin,
\textit{General two-locus selection models: Some objectives,
results and interpretation},
Theor.\ Pop.\ Biol.\ {\bf 7} (1975) 364--398.

\vspace*{-1mm}
\bibitem{KL}
S.~Karlin and U.~Liberman,
\textit{Central equilibria in multilocus systems. I. Generalized
nonepistatic selection regimes},
Genetics {\bf 91} (1979) 777--798.

\vspace*{-1mm}
\bibitem{Kimura}
M.~Kimura,
\textit{A stochastic model concerning the maintenance
of genetic variability in quantitative characters},
Proc.\ Natl.\ Acad.\ Sci.\ {\bf 54} (1965) 731--736.

\vspace*{-1mm}
\bibitem{Kingman}
J.~F.~C.~Kingman,
{\em Markov population processes},
J.\ Appl.\ Prob.\ {\bf 6} (1969) 1--18.

\vspace*{-1mm}
\bibitem{KiLy}
V.~Kirzhner and Yu.~Lyubich,
{\em Multilocus dynamics under haploid selection},
J.\ Math.\ Biol.\ {\bf 35} (1997) 391--408.

\vspace*{-1mm}
\bibitem{Lang}
S.~Lang,
{\em Real and Functional Analysis}, 3rd ed.,
Springer, New York (1993).

\vspace*{-1mm}
\bibitem{Lyu}
Yu.~I.~Lyubich,
\textit{Mathematical Structures in Population Genetics},
Springer, Berlin (1992).

\vspace*{-1mm}
\bibitem{Manos}
H.~Manos and U.~Liberman,
{\em Discrete chiasma formation models and their associated high 
order interference}, J.\ Math.\ Biol.\ {\bf 36} (1998) 448--468.

\vspace*{-1mm}
\bibitem{MHR}
D.~McHale and G.~A.~Ringwood,
{\em Haldane linearisation of baric algebras},
J.\ London Math.\ Soc.~(2) {\bf 28} (1983) 17--26.

\vspace*{-1mm}
\bibitem{RS}
M.~Reed and B.~Simon,
{\em Functional Analysis}, 2nd ed.,
Academic Press, San Diego, CA (1980).

\vspace*{-1mm}
\bibitem{Ringwood}
G.~A.~Ringwood, 
\textit{Hypergeometric algebras and Mendelian genetics},
Niew Archief voor Wiskunde (4) {\bf 3} (1985) 69--83.

\vspace*{-1mm}
\bibitem{Rudin}
W.~Rudin,
{\em Real and Complex Analysis}, 3rd ed.,
McGraw-Hill, New York (1987).

\vspace*{-1mm}
\bibitem{Sch}
H.~H.~Schaefer,
{\em Banach Lattices and Positive Operators},
Springer, Berlin (1974).

\vspace*{-1mm}
\bibitem{Schaeffer}
S.~Schaeffer and E.~L.~Miller,
{\em Estimates of linkage disequilibrium and the recombination
parameter determined from segregating nucleotide sites in the
alcohol dehydrogenase region of Drosophila pseudoobscura},
Genetics {\bf 135} (1993), 541--552.

\vspace*{-1mm}
\bibitem{Colin}
C.~J.~Thompson and J.~L.~McBride,
{\em On {E}igen's theory of the self-organization of matter
and the evolution of biological macromolecules},
Math.\ Biosci.\ {\bf 21} (1974) 127--142.

\vspace*{-1mm}
\bibitem{Werner}
D.~Werner,
{\em Funktionalanalysis}, 3rd ed., Springer, Berlin (2000).

\end{thebibliography}
\end{document}